\documentclass[12pt,letterpaper]{amsart}
\usepackage[margin=1in]{geometry}

\usepackage{amsmath,amssymb,amsthm,graphicx}
\usepackage{xcolor}
\definecolor{referenceblue}{RGB}{0,56,105}
\usepackage[colorlinks=true,linkcolor=referenceblue,citecolor=referenceblue,
            urlcolor=referenceblue,filecolor=referenceblue]{hyperref}
\hypersetup{pdftitle={Effective dynamics of gray solitons in a slowly varying momentum field},pdfauthor={Justin Holmer and Numann Malik}}
\newtheorem{theorem}{Theorem}[section]
\newtheorem{proposition}[theorem]{Proposition}

\newtheorem{corollary}[theorem]{Corollary}
\theoremstyle{remark}

\numberwithin{equation}{section}
\DeclareMathOperator{\sech}{sech}
\newcommand{\X}{\mathcal X}
\newcommand{\PP}{\mathcal P}
\newcommand{\HH}{\mathcal H}
\newcommand{\R}{\mathbb R}

\allowdisplaybreaks[2]
\title[Effective dynamics of gray solitons]{Effective dynamics of gray solitons in a slowly varying momentum field}
\author{Justin Holmer}
\address{Department of Mathematics, Brown University, Providence, RI 02912}
\email{justin\_holmer@brown.edu}
\author{Numann Malik}
\address{Department of Mathematics and Statistics, The College of New Jersey,
Ewing, NJ 08628}
\email{malikn@tcnj.edu}
\date{November 25, 2026}
\begin{document}
\begin{abstract}
We study gray solitons of the one-dimensional defocusing Gross--Pitaevskii equation under a prescribed field $W(hx)$ coupled to the renormalized momentum density. On a compact gray branch we justify the Hamiltonian equations for the center and depth up to times $\kappa h^{-1}\log(1/h)$. The error about the fitted soliton is $O(h^{1/2-\epsilon})$ in the hydrodynamic energy norm, and $O(\sqrt h)$ on bounded slow-time intervals. For exact soliton data we identify the radiation on the scale $z=hx$, $s=ht$. The rescaled remainder converges weakly in $L^2_z$ to two acoustic fields with velocities $W(z)\pm2$, emitted as the core's density deficit and phase jump change. Whenever the depth varies on a fixed slow interval, each hydrodynamic component has a terminal $L^2$ defect of order $\sqrt h$ at distances of order $h^{-1}$ from the core, even if the depth has returned to its initial value. Two explicit fields illustrate the reduced phase portrait and a reversal of the fitted laboratory velocity.

\end{abstract}
\maketitle

\section{Introduction}
In the one-dimensional defocusing Gross--Pitaevskii equation
\begin{equation}\label{intro:gp}
 iu_t=-u_{xx}-2(1-|u|^2)u,\qquad |u(x,t)|\longrightarrow1 \quad(|x|\longrightarrow\infty),
\end{equation}
a dark soliton is a localized depression in the density of the background. Tsuzuki's solitary-wave family \cite{Tsuzuki1971} links its depth, width, intrinsic velocity and phase difference through a single parameter. With
\[
 q_\beta(y)=\cos\beta\,\tanh(y\cos\beta)+i\sin\beta,\qquad
 c(\beta)=2\sin\beta,\qquad -\pi/2<\beta<\pi/2 ,
\]
the function $q_\beta(x-a_0-c(\beta)t)$ solves \eqref{intro:gp}. Its minimum modulus is $|\sin\beta|$. The stationary black soliton $q_0$ has a zero at the center. For $0<\beta<\pi/2$ the phase difference from left to right is $2\beta-\pi$. As $\beta\uparrow\pi/2$, the density deficit vanishes and the intrinsic speed tends to the sound speed $2$.

\begin{figure}[htbp]
\centering
\includegraphics[height=1.95in,width=.43\linewidth,keepaspectratio]{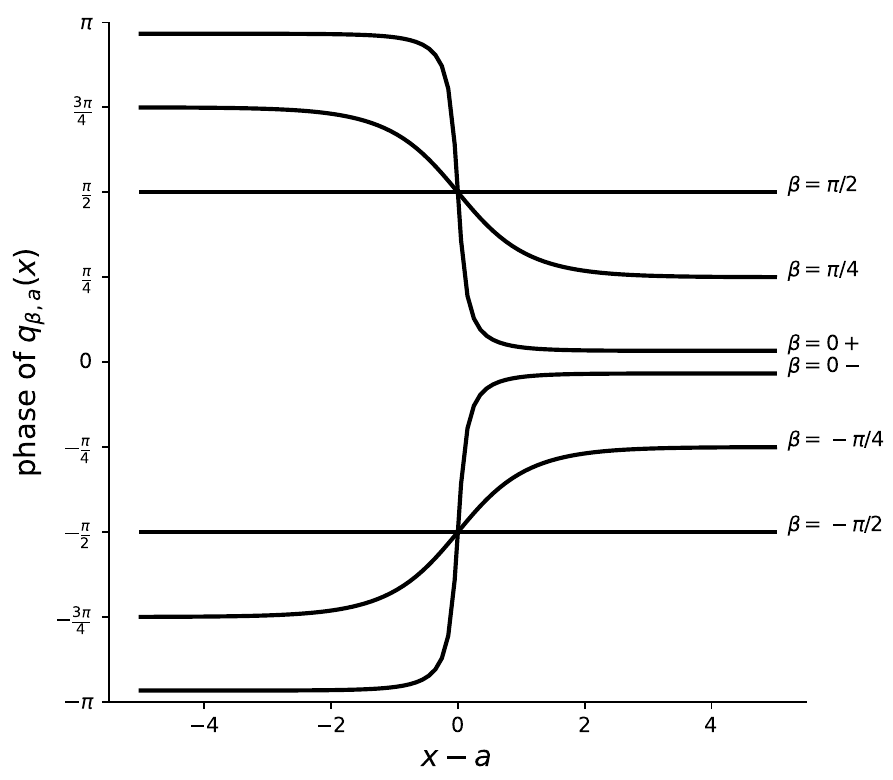}\hspace{.28in}
\includegraphics[height=1.95in,width=.43\linewidth,keepaspectratio]{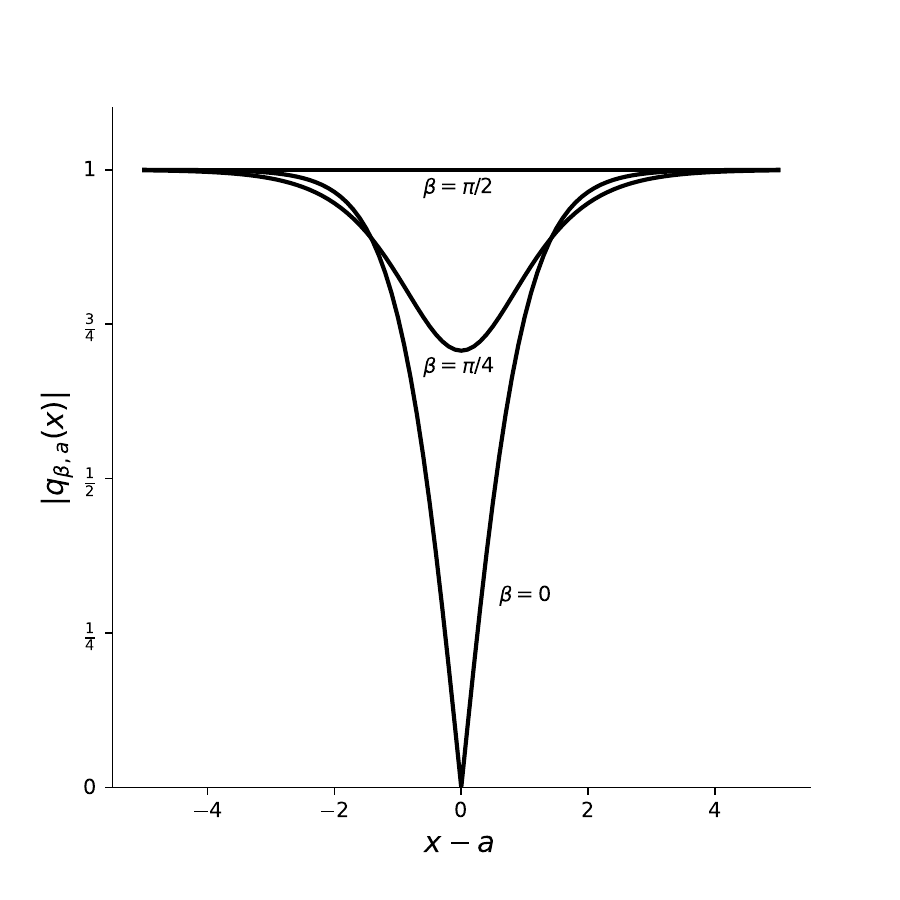}
\caption{Phase and modulus of $q_\beta$. The one-sided limits of the phase at $\beta=0$ record a jump across the zero of the black soliton.}\label{fig:complex}
\end{figure}

We ask how a field varying on length scale $h^{-1}$ moves a gray soliton over a distance of that order and changes its depth by order one. The choice of the field matters because of the nonzero background. In this paper the field couples to the renormalized momentum density. Its value transports the coherent structure, and its gradient changes the momentum of the core. This permits order-one changes in the soliton parameters with a small error in the energy norm, even when the field amplitude is independent of $h$.

\subsection{Hydrodynamical coordinates and the perturbation}
For a nonvanishing field $u=|u|e^{i\phi}$ write
\[
 U=(\eta,\nu)^T,\qquad \eta=1-|u|^2,\qquad \nu=-\phi_x.
\]
We work in the open set $\sup\eta<1$ of $\X=H^1(\R)\times L^2(\R)$, with $\|r\|_{\X}^2=\|r_1\|_{H^1}^2+\|r_2\|_2^2$. The soliton profile is
\begin{equation}\label{intro:Q}
 Q_\beta(y)=
 \begin{bmatrix}\hat\eta_\beta(y)\\ \hat\nu_\beta(y)\end{bmatrix}
 =
 \begin{bmatrix}
 \cos^2\beta\,\sech^2(y\cos\beta)\\[2pt]
 \sin\beta\,\hat\eta_\beta(y)/(1-\hat\eta_\beta(y))
 \end{bmatrix}.
\end{equation}
Both components are localized on a gray branch, and the ambient constant background corresponds to $U=0$.
\begin{figure}[htbp]
\centering
\includegraphics[height=1.98in,width=.43\linewidth,keepaspectratio]{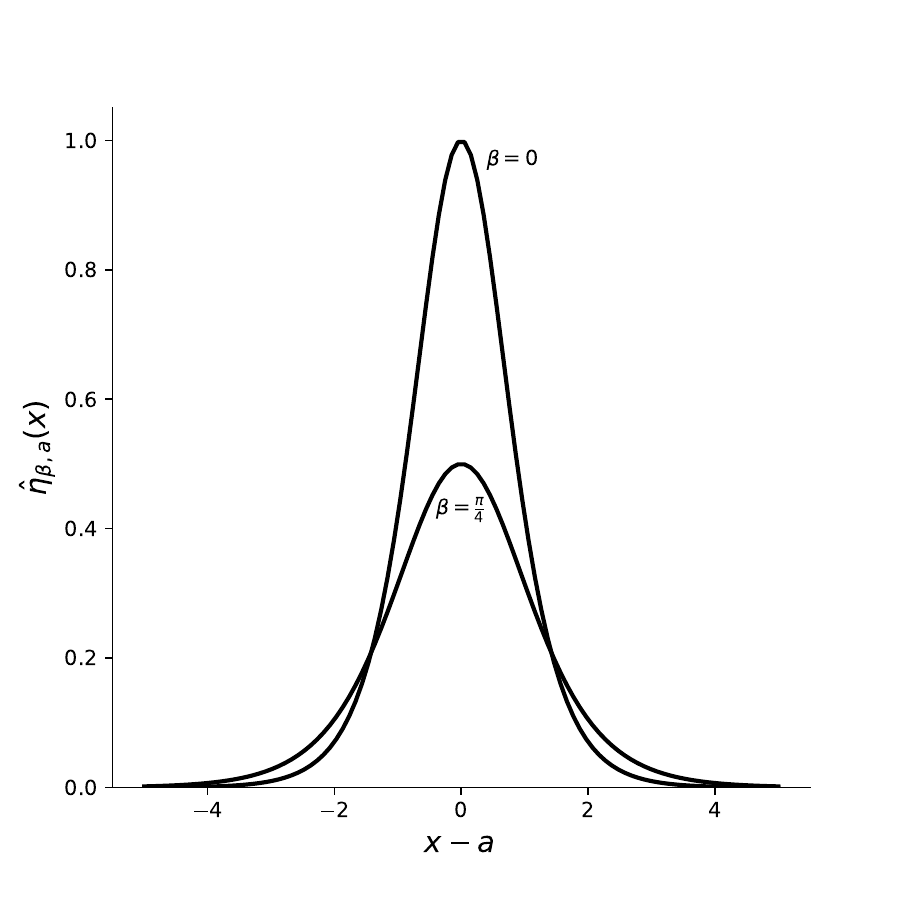}\hspace{.28in}
\includegraphics[height=1.98in,width=.43\linewidth,keepaspectratio]{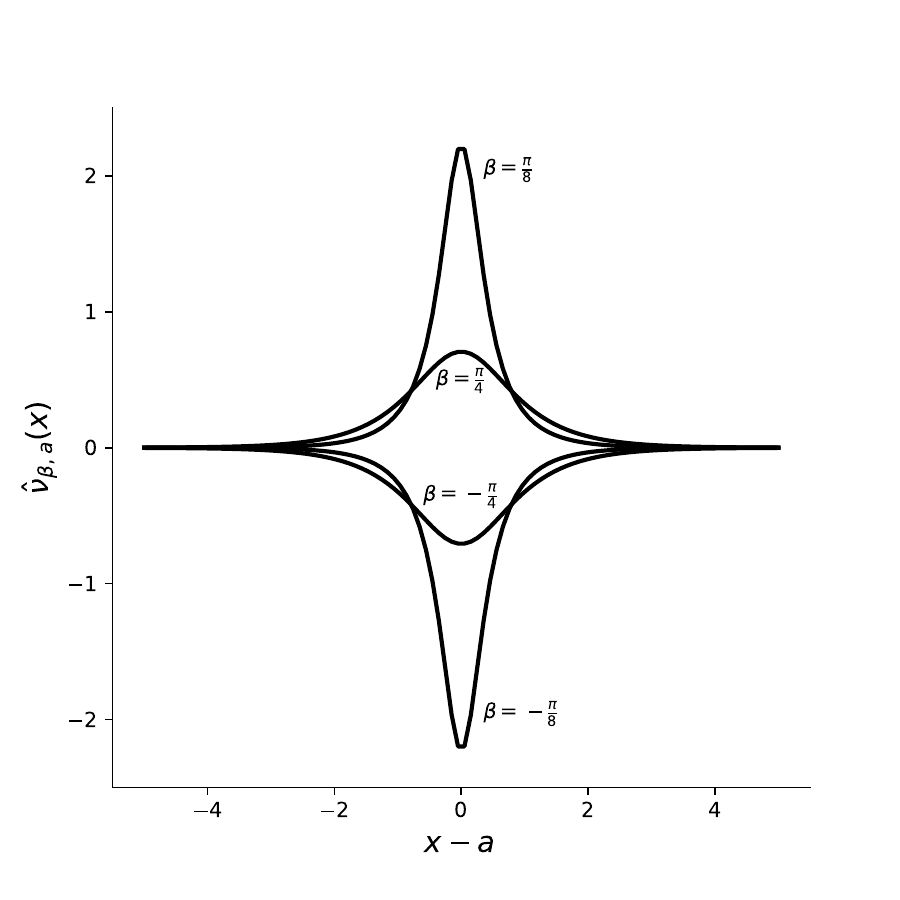}
\caption{Density deficit and negative phase gradient. The density has a smooth black limit, whereas $\hat\nu_\beta(0)=\cos^2\beta/\sin\beta$ diverges.}\label{fig:hydro}
\end{figure}
Let $\sigma=\left(\begin{smallmatrix}0&1\\1&0\end{smallmatrix}\right)$ and $J=-2\partial_x\sigma$. The unforced equation takes the form $U_t=JH_0'(U)$ with energy and renormalized momentum
\begin{equation}\label{intro:energy}
 H_0(U)=\int\left\{\frac{\eta_x^2}{8(1-\eta)}
 +\frac12(1-\eta)\nu^2+\frac12\eta^2\right\}dx,\qquad
 \PP(U)=\frac12\int\eta\nu\,dx .
\end{equation}
All spatial integrals without a domain are over $\R$. Fix a real $W\in\mathcal S(\R)$ and take $V(x)=W(hx)$ with $0<h\leq1$. Our perturbed Hamiltonian and flow are
\begin{equation}\label{intro:forced}
 H(U)=H_0(U)+\tfrac12\int V\eta\nu\,dx,\qquad
 U_t=JH_0'(U)-\partial_x(VU).
\end{equation}
Since $J\PP'(U)=-U_x$, freezing $V$ at the center produces an exact translation. In centered coordinates $y=x-a$, the rest of the field acting on $Q_\beta$ is the localized source $-\partial_y[(V(a+y)-V(a))Q_\beta]$, of size $O(h)$. The constant background remains an exact solution. By comparison, the scalar density coupling $\tfrac12\int V\eta$ drives $\nu$ by $-V_x$ throughout the region where the field varies. The coupling $\tfrac12\int V\eta^2$ preserves the background, but its frozen vector field at the soliton is $(0,-2V(a)\hat\eta_\beta')^T$. For $V(a)\ne0$ this vector lies outside the tangent space. Indeed its odd parity rules out a depth component, and the translation tangent has nonzero first component.
In complex variables \eqref{intro:forced} reads
\begin{equation}\label{intro:complex}
 iu_t=-u_{xx}-2(1-|u|^2)u-iVu_x-\tfrac i2 V_x(u-1/\bar u).
\end{equation}
This expression also records an obstruction at a zero of $u$ where $V_x\ne0$.

\subsection{The effective equations}
On the positive gray branch write $E(\beta)=H_0(Q_\beta)$ and $P(\beta)=\PP(Q_\beta)$. Direct integration gives
\begin{equation}\label{intro:EP}
 E(\beta)=\tfrac43\cos^3\beta,\quad
 P(\beta)=\tfrac\pi2-\beta-\sin\beta\cos\beta,\quad
 P'(\beta)=-2\cos^2\beta,\quad \frac{E'(\beta)}{P'(\beta)}=2\sin\beta .
\end{equation}
Indeed the energy integrand in complex variables equals $\cos^4\beta\,\sech^4(y\cos\beta)$ and $\int\operatorname{Im}(\bar q_\beta q_\beta')=-2\sin\beta\cos\beta$. Subtracting the phase difference from the latter integral and dividing by two gives $P$. For an arbitrary energy datum, \eqref{intro:energy} defines the momentum even when separate phase limits are unavailable.

Put $p_\beta(y)=\hat\eta_\beta(y)\hat\nu_\beta(y)/2$ and
$M_2(\beta)=\int y^2p_\beta(y)\,dy$. On the positive gray branch the restricted Hamiltonian has the expansion
\begin{equation}\label{eq:finite-width-hamiltonian}
 H(Q_\beta(\cdot-a))
 =E(\beta)+W(ha)P(\beta)
  +\tfrac12 h^2 W''(ha)M_2(\beta)+\mathcal R_h(a,\beta).
\end{equation}
On each compact subinterval $I$ of this branch,
\[
 |\partial_\beta^j\partial_a^\ell\mathcal R_h(a,\beta)|
 \le C_{I,W}h^{4+\ell},\qquad j,\ell\in\{0,1\},\quad 0<h\le1 .
\]
Indeed, the field energy equals $\int W(ha+hy)p_\beta(y)\,dy$.
The first and third moments of $p_\beta$ and $\partial_\beta p_\beta$
vanish by parity, and their absolute fourth moments are bounded uniformly
on $I$. Taylor expansion through degree three gives the estimate.
An $a$ derivative replaces $W$ by $hW'$ in the same argument.

Momentum generates translations, since $J\PP'(Q_\beta)=\partial_a Q_\beta(\cdot-a)$ in translated coordinates. On the exponentially localized tangents we use $\partial_y^{-1}f(y)=\tfrac12\int\operatorname{sgn}(y-z)f(z)\,dz$, which is skew under the integral pairing. Consequently the restricted two-form is $dP\wedge da$ in the convention $\omega(F,G)=\langle F,J^{-1}G\rangle$. Its coefficient is determined by the localized pairing $\langle\partial_\beta Q_\beta,\sigma Q_\beta/2\rangle=P'(\beta)$; no bounded inverse of $J$ on $\X$ is required. In slow time $s=ht$ with $A$ the scaled center, the effective Hamiltonian is $\HH(A,B)=E(B)+W(A)P(B)$. Its equations are
\begin{equation}\label{ode}
 A'=2\sin B+W(A),\qquad
 \frac{d}{ds}P(B)=-W'(A)P(B),\qquad
 B'=\frac{W'(A)P(B)}{2\cos^2 B}.
\end{equation}
Here $A'$ is the velocity measured in the laboratory. Since $A=ha$ and $s=ht$ have the same rescaling factor, it is compared with $\dot a$ in physical time. The proof below also recovers the force from the exact PDE momentum balance.

Ivanov and Kamchatnov derived a Hamiltonian of this form for motion through a nonuniform condensate \cite[equations~(14)--(17)]{IvanovKamchatnov2022}. At unit density and relative speed $\sin B$, their canonical momentum is $p=2P(B)$ and their intrinsic energy is $E(B)$. Taking their background velocity $u=W/2$ gives $up=WP$ in our normalization. In their model the density and velocity evolve as background fields. In \eqref{intro:forced} the coefficient $W$ is prescribed and the unit background is an exact solution. Theorem~\ref{thm:main} quantifies the accuracy of \eqref{ode} for this coupling.

\subsection{The finite-time theorem}
Write $r(y,t)=U(a(t)+y,t)-Q_{\beta(t)}(y)=(\hat\theta,\hat\xi)^T$. Fit the parameters by
\begin{equation}\label{intro:orth}
 \langle r,Q_\beta'\rangle=\langle r,\sigma Q_\beta\rangle=0
\end{equation}
with spatial derivative $Q_\beta'$ and the real $L^2$ pairing. These conditions are well defined throughout a small $\X$-neighborhood of the gray family on a compact parameter interval.

\begin{theorem}[Effective dynamics]\label{thm:main}
Fix a real Schwartz field $W$, closed intervals $K\subset\operatorname{int}I$, $I\Subset(0,\pi/2)$, and $C_0>0$. Let $0<\epsilon<1/2$. Suppose the initial datum $U_0=(\eta_0,\nu_0)^T\in\X$ satisfies $\inf_x(1-\eta_0)>0$ and has a fitted decomposition \eqref{intro:orth} with $\beta(0)\in\operatorname{int}K$ and $\|r(0)\|_{\X}\leq C_0\sqrt h$. Solve \eqref{ode} from $(A(0),B(0))=(ha(0),\beta(0))$ and define
\[
 S_K=\inf\{s\geq0:B(s)\notin\operatorname{int}K\},\qquad
 T_h=h^{-1}\min\{S_K,\kappa\log(1/h)\},\qquad\inf\varnothing=\infty .
\]
There are constants $C,R,\kappa,h_0>0$, independent of the initial center, such that for $0<h<h_0$ the energy solution exists and is nonvanishing on $[0,T_h]$, with continuously differentiable fitted parameters and $\beta(t)\in I$. Write
\[
 D_R(t)=\int \sech^2(y/R)
   \bigl(\hat\theta_y^2+\hat\theta^2+\hat\xi^2\bigr)(y,t)\,dy .
\]
For every $T\leq T_h$,
\begin{align}
 \sup_{0\leq t\leq T}\|r(t)\|_{\X}^2+\int_0^T D_R(t)\,dt
 &\leq C(\|r(0)\|_{\X}^2+h^2 T)e^{ChT},\label{intro:spacetime}\\
 |ha(T)-A(hT)|+|\beta(T)-B(hT)|
 &\leq Ch e^{ChT}.\label{intro:trajectory}
\end{align}
In particular $\sup_{[0,T_h]}\|r\|_{\X}\leq C h^{1/2-\epsilon}$ with $\kappa$ chosen in terms of $\epsilon$. The field error is $O(\sqrt h)$ on every fixed slow interval covered by the theorem. For exact soliton data the right side of \eqref{intro:trajectory} can be sharpened to $C h^2 T e^{ChT}$.
\end{theorem}

An $O(\sqrt h)$ neighborhood of a reference soliton gives fitted initial parameters by modulation, with an $O(\sqrt h)$ adjustment. This adjustment matters for the trajectory estimate. For example, if $W=0$, an exact soliton with angle $\beta_0+\sqrt h$ and the soliton of angle $\beta_0$ have centers separated by order $h^{-1/2}$ at time $1/h$ when they start at the same point. The field estimate \eqref{intro:spacetime} uses the profile centered at $a(t)$. On a bounded slow interval, \eqref{intro:trajectory} controls $ha$ to order $h$, allowing an order-one error in the physical center.

There is a second conclusion on fixed slow intervals. For the exact initial soliton with scaled center $A_0$ and angle $B_0$, let $R_h(s,z)$ denote $h^{-1}$ times the laboratory-frame remainder at $(t,x)=(s/h,z/h)$. Theorem~\ref{thm:acoustic} identifies its macroscopic limit. The two scalar components $R^{\pm}=R_1\pm R_2$ solve
\[
 \partial_s R^{\pm}+\partial_z\big((W(z)\pm2)R^{\pm}\big)
     =-\frac{d}{ds}\big[2\cos B(s)\pm(\pi-2B(s))\big],\delta_{A(s)},
 \qquad R^{\pm}(0)=0.
\]
The quantities in brackets are sums and differences of the integrated density deficit and the integrated negative phase gradient of the core. Emissions at different times along either family of characteristics arrive at different positions. This yields the terminal lower bounds of Corollary~\ref{cor:acoustic-lower} even when $B(S)=B(0)$.

If the effective trajectory remains in a compact gray interval, one can choose $K$ to contain it in its interior, obtaining the full logarithmic comparison. If the angle tends to zero at a finite slow time the result extends to each fixed compact gray portion preceding that event. The constants depend on the gray margin. For a fixed nonzero energy level $e=\HH(A,B)$ the sonic boundary is uniformly separated, since $0\leq E(B)\leq 2P(B)$ by \eqref{intro:EP} and thus $P(B)\ge |e|/(2+\|W\|_\infty)$. Reflection in space, with $W(A)$ replaced by $-W(-A)$, gives the negative-branch theorem. Intrinsic speed tending to zero and laboratory speed tending to zero have distinct dynamical consequences, as the next examples illustrate.

\subsection{Two examples in the reduced phase plane}
Consider a soliton incident from the left on the positive gray branch, with limiting angle $B_{\rm in}\in(0,\pi/2)$ as $A\to-\infty$. Its energy level is $e=E(B_{\rm in})\in(0,4/3)$. The incoming and outgoing states here refer to the orbits of \eqref{ode}. Theorem~\ref{thm:main} compares a PDE solution started at finite slow time with the corresponding portion of the orbit.
\begin{figure}[tbp]
\centering
\includegraphics[width=.84\linewidth,height=5.70in,keepaspectratio]{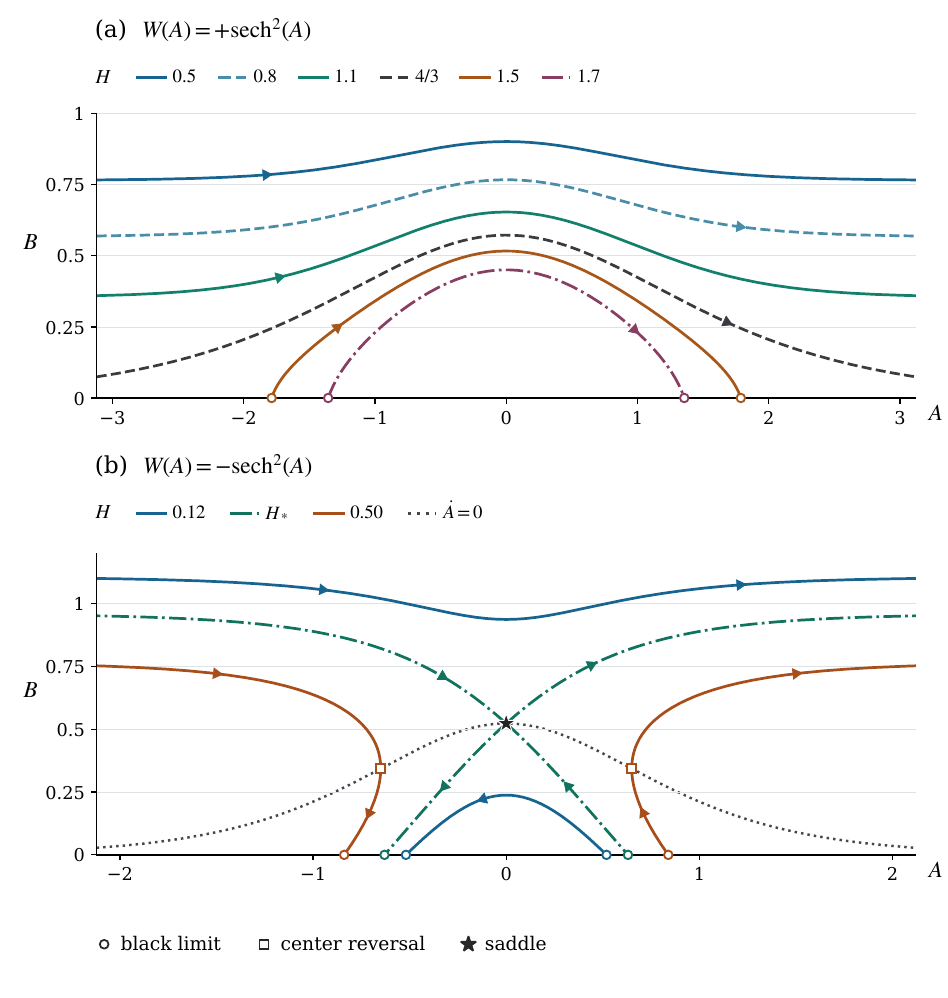}
\caption{The Hamiltonian level curves of \eqref{ode}, with arrows in increasing slow time. In the upper panel incident gray orbits cross the positive field. In the lower panel the dotted curve gives zero laboratory speed. The incoming orbits with energies $0.12,e_*,0.50$ transmit, approach the saddle, and turn back, respectively. A square marks reversal at positive $B$ and an open circle marks intersection with the black boundary.}\label{fig:phaseplanes}
\end{figure}

\emph{Aiding field $W(A)=+\sech^2 A$.}
Since $E$ decreases and $e=E(B)+W(A)P(B)\ge E(B)$, one has $B\ge B_{\rm in}$. In particular, $A'\ge2\sin B_{\rm in}>0$. As the orbit approaches the center from the left $B$ increases. After $A=0$ it decreases and returns to $B_{\rm in}$ as $A\to+\infty$. The whole trajectory has a uniform gray margin. At fixed $A$ the energy is strictly decreasing in $B$ from $4/3+(\pi/2)W(A)$ to zero. The level $e=4/3$ therefore tends to $B=0$ at both spatial infinities, and levels $4/3<e<4/3+\pi/2$ give arcs with black endpoints at finite positions.

\emph{Opposing field $W(A)=-\sech^2 A$.}
Set $v=\sech^2 A$ and $K_v(B)=E(B)-vP(B)$. Since $K_v'(B)=-2\cos^2 B(2\sin B-v)$, its unique interior maximum for $v>0$ is
\[
 M(v)=K_v(\arcsin(v/2)),\qquad M'(v)=-P(\arcsin(v/2))<0 .
\]
The least maximum occurs at the field center. The critical incoming energy is
\[
 e_*=M(1)=\frac{3\sqrt3}{4}-\frac{\pi}{3}=0.25184055\ldots .
\]
For $0<e<e_*$ the incoming orbit lies on the root $2\sin B>v$ throughout. Continuity of this root for $0\le v\le1$ gives a uniform gray margin and a strictly positive lower bound for $A'$. The orbit transmits and recovers $B_{\rm in}$.

For $e=e_*$ the incoming root tends to the equilibrium $(A,B)=(0,\pi/6)$. There the mixed energy derivative vanishes, $\HH_{AA}=W''(0)P(\pi/6)>0$ and $\HH_{BB}=-4\cos^3(\pi/6)<0$. Thus this is a hyperbolic equilibrium. The incident orbit approaches it in infinite slow time by uniqueness of the smooth ODE. Its laboratory speed tends to zero while its minimum amplitude tends to $1/2$. The full logarithmic comparison applies.

For $e_*<e<4/3$ the center reaches the unique left-flank point with $M(v_{\rm turn})=e$. It reaches this point in finite slow time since $B'$ is strictly negative and bounded away from zero near it. Here $B_{\rm turn}=\arcsin(v_{\rm turn}/2)>0$, $A'=0$, and $A''=2\cos B_{\rm turn}\,B'<0$. The orbit then follows the other root, with $A'<0$ and $B'<0$, to the black limit at
\[
 \sech^2 A_{\rm black}=\frac{2}{\pi}\Big(\frac43-e\Big),\qquad
 A_{\rm black}<A_{\rm turn}<0 .
\]
Indeed $M(v_{\rm turn})>K_{v_{\rm turn}}(0)$ gives this ordering. Along this root the limiting value of $B'$ at $A_{\rm black}$ is strictly negative, giving finite slow time there as well. Corollary~\ref{cor:reversal} gives a reversal of the fitted PDE velocity on each fixed interval containing this turn in its interior and remaining in the gray band. The approach to zero eventually leaves every fixed gray band. Transmission, trapping and reversal also occur for dark solitons incident on evolving hydrodynamic waves in defocusing NLS, as studied by Sprenger, Hoefer and El \cite{SprengerHoeferEl2018}.

\subsection{Relation to previous work}
Kivshar and Yang derived adiabatic equations for the phase angle of
perturbed dark solitons on constant and varying backgrounds
\cite{KivsharYang1994}. For trapped condensates, Busch and Anglin used
matched asymptotics to obtain the motion of a dark soliton through a
slowly varying background density and flow \cite{BuschAnglin2000}.
Konotop and Pitaevskii treated gray solitons of arbitrary amplitude as
quasiparticles in a local-density approximation
\cite{KonotopPitaevskii2004}. Pelinovsky, Frantzeskakis and Kevrekidis
derived radiative corrections to the oscillation law in a parabolic
trap \cite{PelinovskyFrantzeskakisKevrekidis2005}. For a small decaying
scalar potential, Pelinovsky and Kevrekidis proved persistence results
for pinned black solitons, characterized their instabilities near
extrema of an effective potential and compared their numerical evolution
with a particle law
\cite{PelinovskyKevrekidis2008}.

For the prescribed coupling to the renormalized momentum density in
\eqref{intro:forced}, we justify the Hamiltonian approximation
\eqref{ode} with bounds on the macroscopic center, the angle and the
hydrodynamic energy norm. The comparison extends to the first exit
of the reduced angle from the chosen gray interval or the logarithmic
time bound. For exact soliton data we also determine the weak
macroscopic limit of the density and phase-gradient remainder.
We relate this result to earlier work on radiation shelves in Section~\ref{sec:acoustic}.

In the unforced equation, B\'ethuel, Gravejat and Smets proved
orbital stability for well-separated chains of gray solitons and
asymptotic stability for a single gray soliton
\cite{BGSChains2014,BGS_Asymp_Dark_2015}. Gravejat and Smets established
asymptotic stability for the black soliton \cite{GS_Asymp_Black_2015}.
At the black limit the hydrodynamic chart degenerates, and the forced
complex equation \eqref{intro:complex} contains a reciprocal of the field.

Fr\"ohlich, Gustafson, Jonsson and Sigal justified effective particle
dynamics for localized nonlinear Schr\"odinger solitary waves in an
external potential \cite{FGJS2004}. Holmer and Zworski obtained a
logarithmic-time Hamiltonian approximation for a bright soliton in a
slowly varying scalar potential \cite{HolmerZworski2008}. Our use of a
local virial and a modulated energy follows Holmer's treatment of
slowly varying transport in KdV on the $h^{-1}\log(1/h)$ time scale
\cite{Holmer2011}, building on the local virial method of Martel and
Merle \cite{MartelMerle2005}.

For the unforced Gross--Pitaevskii equation, B\'ethuel, Gravejat and
Smets applied a Hessian transform and a virial estimate to a limiting
orbit \cite[Section~1.3.3 and Proposition~7]{BGS_Asymp_Dark_2015}.
Virial estimates on both the perturbation and a regularized transform
also enter the nonlinear Klein--Gordon and cubic--quintic
Schr\"odinger arguments of \cite{KMM2022,MartelCQ2022}. Our companion
\cite{HM} supplies the spatial estimates for the gray profile, including
coercivity of the dual commutator with weight $y(1-\hat\eta_\beta)$.
Here the depth is driven at order $h$, and the perturbation can reach
points where the field differs by order one from its value at the core.
We localize the transport term before applying the dual transformation
and use the exact total momentum balance to estimate the parameters.

Section~\ref{sec:flow} constructs the forced nonvanishing flow. Sections~\ref{sec:fitted-modulation}--\ref{vir:section} establish the modulation, action and virial estimates used in Section~\ref{sec:closure} to prove Theorem~\ref{thm:main}. Section~\ref{sec:sharpness} proves a terminal lower bound directly from the fluxes when the endpoint depths differ. Section~\ref{sec:acoustic} identifies the limiting radiation, including its dependence on any prepared initial field.

\subsection{Disclosure on use of AI}
The authors chose the perturbation coupling the external field to
the renormalized momentum density, derived the effective equations
by symplectic restriction, and worked out the phase portraits for
$W(A)=\pm\sech^2 A$. Following \cite{Holmer2011}, we chose to combine
a Lyapunov functional with local virial estimates. Our preliminary
calculations indicated that this approach should give an $O(\sqrt h)$
bound for the remainder in the hydrodynamic energy norm on bounded
slow-time intervals.  We used ChatGPT Astra 6.2 to help complete 
the energy and virial
estimates, review the proofs and their relation to earlier work, and
draft and revise the exposition. During revision, 
ChatGPT formulated and proved Theorem~\ref{thm:acoustic} 
and the resulting terminal
lower bounds in Corollary~\ref{cor:acoustic-lower}.

\section{The perturbed flow and its balance laws}\label{sec:flow}
Put $\mathcal X=H^1(\mathbb R)\times L^2(\mathbb R)$ and
$\mathcal O=\{(\eta,\nu)\in\mathcal X:\inf_x(1-\eta(x))>0\}$.
For this section $h>0$ is fixed, so $V(x)=W(hx)$ is a fixed real Schwartz function.
Write $\rho=\sqrt{1-\eta}$. The Hamiltonian and momentum are
\begin{equation}\label{eq:flow-ham}
 H(U)=\int\left(\frac{\eta_x^2}{8\rho^2}
       +\frac12\rho^2\nu^2+\frac12\eta^2+\frac12 V\eta\nu\right)dx,
 \qquad \mathcal P(U)=\frac12\int\eta\nu\,dx .
\end{equation}
With $J=-2\partial_x\sigma$, where $\sigma$ exchanges the two components,
the equations $U_t=JH'(U)$ read
\begin{equation}\label{eq:flow-hydro}
 \eta_t=-\partial_x(2\rho^2\nu+V\eta),\qquad
 \nu_t=-\partial_x G,\qquad
 G=\frac{\rho_{xx}}{\rho}-\nu^2+2\eta+V\nu .
\end{equation}
Indeed $H'_\nu=\rho^2\nu+V\eta/2$ and $2H'_\eta=G$.
Here $\rho-1\in H^1$ and $G\in H^{-1}$, using the $H^1$ multiplier action
on $H^{-1}$ in dimension one for the nonconstant part of $1/\rho$.

For $u=\rho e^{i\phi}$ and $\nu=-\phi_x$, the complex equation is
\begin{equation}\label{eq:flow-complex}
 i u_t=-u_{xx}+2(|u|^2-1)u-iVu_x
             -\frac{i}{2}V_x\left(u-\frac{1}{\bar u}\right).
\end{equation}
To see the last coefficient, write the added Hamiltonian density as
$V(\operatorname{Im}(\bar u u_x)-\phi_x)/2$.
For a compactly supported variation $v$, the phase variation is
$\operatorname{Im}(v/u)$. Integrating by parts gives the real $L^2$ gradient
$-iVu_x-iV_x u/2+iV_x/(2\bar u)$.
Taking the real and imaginary parts of \eqref{eq:flow-complex} in polar
coordinates gives \eqref{eq:flow-hydro} and $\phi_t=G$. For a constant
field $V_*$, the solution is the unforced solution evaluated at
$(x-V_*t,t)$.

For the unforced equation, Mohamad proved local well-posedness up to
vacuum in the hydrodynamic energy space and its integer-order analogues
\cite[Theorem~1.2]{Mohamad2014}. Wegner established local bi-Lipschitz
estimates for the Madelung transform and global hydrodynamic well-posedness
under energy assumptions excluding vacuum
\cite[Theorems~1.6 and~1.9]{Wegner2023}. We construct the flow for the
momentum field in \eqref{eq:flow-ham} and derive its balance laws below.

\begin{proposition}[Nonvanishing energy flow]\label{prop:flow}
For each $U_0\in\mathcal O$ there is a unique maximal distributional
solution $U\in C((-T_-,T_+);\mathcal O)$ of \eqref{eq:flow-hydro}
with $U(0)=U_0$. The flow is continuous on each compact subinterval of
existence, and $H^m\times H^{m-1}$ regularity persists there for integers
$m\geq1$. If $T_+<\infty$, then
\[
 \lim_{t\uparrow T_+}\inf_{x\in\mathbb R}(1-\eta(t,x))=0.
\]
The same limit holds as $t\downarrow -T_-$ if $T_-<\infty$. Smooth data approximate the
solution in $C_t\mathcal X$ on each compact subinterval containing zero.
On the interval of existence $H(U(t))=H(U_0)$, and, for each
$\chi\in C_c^\infty(\mathbb R)$,
\begin{align}
 \frac{d}{dt}\int\chi\eta\,dx
   &=\int\chi'(2(1-\eta)\nu+V\eta)\,dx,\label{eq:density}\\
 \frac{d}{dt}\mathcal P(U(t))
   &=-\frac12\int V_x\eta\nu\,dx.\label{eq:momentum}
\end{align}
\end{proposition}

\begin{proof}
\emph{Local construction.}
Lift the initial data by taking $u_0=\rho_0 e^{i\phi_0}$,
$\phi_0(x)=-\int_0^x\nu_0(z)\,dz$. Then $u_0$ is bounded,
$u_{0,x}\in L^2$, and $1-|u_0|^2\in L^2$. Set
$u=e^{i\psi}w$ with $\psi(0)=0$ and $\psi'=-V/2$.
The gauge preserves these three properties, and \eqref{eq:flow-complex}
becomes
\begin{equation}\label{eq:flow-gauge}
 iw_t=-w_{xx}+F(x,w),\qquad
 F(x,z)=2(|z|^2-1)z-\tfrac14 V(x)^2z+
                             \frac{i V_x(x)}{2\bar z}.
\end{equation}
Let $S(t)=e^{it\partial_x^2}$. For any bounded $f$ with $f_x\in L^2$,
$d_f(t)=(S(t)-1)f\in C_tH^1$, with
\[
 \|d_f(t)\|_2\leq C|t|^{1/2}\|f_x\|_2,\qquad
 \|(d_f)_x(t)\|_2\leq2\|f_x\|_2,\qquad
 \|d_f(t)\|_\infty\leq C|t|^{1/4}\|f_x\|_2 .
\]
The first bound follows by applying the multiplier
$(e^{-itk^2}-1)/(ik)$ to $f_x$, and the third follows by interpolation.
For $f=e^{-i\psi}u_0$, look for $w=f+d_f+z$ with
\[
 z(t)=-i\int_0^t S(t-\tau)F(\cdot,f+d_f(\tau)+z(\tau))\,d\tau .
\]
On a bounded $C_tH^1$ ball for $z$, chosen small in $L^\infty_{t,x}$
relative to $\inf|f|$, the argument of $F$ stays in a fixed annulus
on a sufficiently short interval. On this ball $F$ is bounded in
$C_t H^1$ and Lipschitz with respect to $z$ in that norm. Indeed,
$|w|^2-1\in L^2$, $w_x\in L^2$, and $V^2,V_x\in H^1$. The first
derivatives of the smooth functions of $w$ are bounded on the annulus.
For differences, their coefficient differences are bounded pointwise
by $C|z_1-z_2|$, using a smooth extension of the reciprocal across the
excluded disk. The product rule and $H^1\hookrightarrow L^\infty$ give
the stated Lipschitz bound. The Duhamel map therefore contracts on a time
interval depending only on $V$, an upper bound for $\|U_0\|_{\mathcal X}$
and a positive lower bound for $1-\eta_0$.
Gauge equivariance gives the same hydrodynamic solution for every choice
of constant initial phase. The same local existence time applies on
restarting from any datum with these bounds.

\emph{Continuity in hydrodynamic coordinates.}
Let $U_{0,n}\to U_0$ in $\mathcal X$ and normalize the initial lifts at
the origin as above. For their gauged lifts $f_n,f$ one has
\[
 f_n\to f\ \text{locally uniformly},\qquad
 f_{n,x}\to f_x\ \text{in }L^2,\qquad
 |f_n|^2-1\to |f|^2-1\ \text{in }L^2 .
\]
In fact the phase differences tend uniformly to zero on compact sets by
Cauchy--Schwarz. In the derivative formula
$f_{n,x}=e^{i(\phi_{0,n}-\psi)}
[\rho_{0,n,x}-i\rho_{0,n}(\nu_{0,n}+\psi')]$, the bracket converges in
$L^2$ and multiplication of its limit by the phase factors converges by
dominated convergence. The free increments $d_{f_n}$ consequently converge
in $C_tH^1$ to $d_f$. To compare the contractions, evaluate both maps
at the fixed solution $z$ for $f$. Put $w_n^*=f_n+d_{f_n}+z$.
These functions converge locally uniformly to $w$, their derivatives
converge in $C_tL^2$, and $|w_n^*|^2-1\to |w|^2-1$ in $C_tL^2$.
For the last assertion, expand the squares about $f_n$ and use the
convergence of $d_{f_n}+z$ in $C_tH^1$. A uniformly bounded coefficient
converging locally uniformly can be multiplied into a fixed
$C_t L^2$ factor with convergence in $C_tL^2$ on a compact time
interval: its time range has uniformly small spatial $L^2$ tails.
Applying this observation to \eqref{eq:flow-gauge} and its spatial
derivative gives $F(\cdot,w_n^*)\to F(\cdot,w)$ in $C_tH^1$.
The uniform contraction gives $z_n\to z$ in $C_tH^1$. Undoing the gauge,
the same product argument applied to
$\eta_x=-2\operatorname{Re}(\bar u u_x)$ and
$\nu=-\operatorname{Im}(u_x/u)$ gives convergence in $C_t\mathcal X$.
Iteration covers each compact subinterval of existence. The argument
allows lifts with different phases at spatial infinity; no integrability
condition on the phase difference enters.

\emph{Distributional correspondence and uniqueness.}
Here are the details at energy regularity. For smooth data the polar
identities follow by differentiation. The complex construction preserves
higher regularity: after $m$ spatial derivatives of \eqref{eq:flow-gauge},
the one-dimensional product and composition estimates bound
$\|F(\cdot,w)\|_{H^m}$ by
$C_m(1+\|w_x\|_{H^{m-1}})$ when the energy bounds, the annulus and $V$
are fixed. This follows from the product rule and the homogeneous
Gagliardo--Nirenberg inequalities for the intermediate derivatives;
the undifferentiated cubic term uses $|w|^2-1\in L^2$.
Duhamel's formula and Gronwall give persistence, which transfers by the
polar formulas to $H^m\times H^{m-1}$. Approximate $U_0$ by smooth
compactly supported hydrodynamic data with a common positive density
floor. The preceding continuity gives convergence in $C_t\mathcal X$.
In particular $G_n\to G$ in $C_tH^{-1}$ since multiplication by
$1/\rho_n$ converges in multiplier norm there and
$\nu_n^2\to\nu^2$ in $C_tL^1$. Thus the limiting pair satisfies
\eqref{eq:flow-hydro} in distributions.

Conversely, given a distributional solution in $C_t\mathcal O$, choose a
spatial primitive $\phi_x=-\nu$ with weighted mean zero against
a fixed smooth compactly supported function of integral one, denoted by
$\omega$. Its time derivative differs from $G$ by the spatially constant
distribution $-\langle G,\omega\rangle$. Adding
$\int_0^t\langle G(s),\omega\rangle\,ds$ and a constant to $\phi$
fixes $\phi_t=G$ and the initial phase. Then
$\rho,\phi\in C_tH^1_{\rm loc}$ with time derivatives in
$C_t H^{-1}_{\rm loc}$. For $\rho$ this follows from
$\rho_t=-\eta_t/(2\rho)$. This chain rule, as well as the time product
rule for $u=\rho e^{i\phi}$, is justified by mollifying in time on compact
subintervals: the arguments converge in $C_t H^1_{\rm loc}$, and their
mollified time derivatives in $C_t H^{-1}_{\rm loc}$, so the products
converge by the $H^1$ multiplier bound after spatial localization.
Substitution gives \eqref{eq:flow-complex} in distributions.
Its gauge $w$ is bounded, has $w_x\in C_tL^2$ by the locally uniform
phase convergence and the tail argument above, and obeys
\eqref{eq:flow-gauge} with $F(\cdot,w)\in C_t H^1$. Variation of
constants in tempered distributions gives the same Duhamel formula,
and its Duhamel correction tends to zero in $H^1$ at the initial time. It therefore lies in the contraction ball for a short time, and iteration gives uniqueness of $U$.

\emph{Balance laws.}
For smooth solutions the skew symmetry of $J$ gives conservation of
the Hamiltonian in \eqref{eq:flow-ham}. The first equation in
\eqref{eq:flow-hydro} gives \eqref{eq:density}.
Since $J\mathcal P'(U)=-U_x$,
\[
 \frac{d}{dt}\mathcal P(U)=\langle H'(U),U_x\rangle
       =\frac12\int V\partial_x(\eta\nu)\,dx .
\]
The unforced term vanishes by translation invariance. An integration by
parts yields \eqref{eq:momentum}. Smooth approximation passes all three
integrated identities to the flow in $\mathcal O$: $H$ and $\mathcal P$
are continuous there, the density flux is continuous into $L^2$, and
$\eta\nu$ is continuous into $L^1$. The right sides of the asserted
differential identities are continuous in time.

\emph{Continuation up to vacuum.}
Suppose that $\inf_x\rho(t,x)^2\geq\mu>0$ at some time $t$. Young's
inequality in the momentum term of \eqref{eq:flow-ham} gives
\[
 H(U(t))\geq \tfrac12\|\rho_x\|_2^2+\tfrac14\|\rho\nu\|_2^2
       +\tfrac12\|\eta\|_2^2-\tfrac14\int V^2\frac{\eta^2}{\rho^2}\,dx .
\]
Since $\eta^2/\rho^2=-1-\eta+\rho^{-2}\leq|\eta|+\mu^{-1}$,
the last integral is bounded by
$\|V^2\|_2\|\eta\|_2+\mu^{-1}\|V\|_2^2$. Absorbing the linear
term in $\|\eta\|_2$ and using conservation of $H$ bounds
$\|\rho_x\|_2$, $\|\rho\nu\|_2$ and $\|\eta\|_2$ in terms of
$H(U_0)$, $\mu$ and $V$. The pointwise inequality
$|\rho-1|\leq|\eta|$ and Sobolev embedding now bound
$\|\rho\|_\infty$. Consequently $\eta_x=-2\rho\rho_x$ and
$\|\nu\|_2\leq\mu^{-1/2}\|\rho\nu\|_2$ give a bound for
$\|U(t)\|_{\mathcal X}$ with the same dependence.
If $t_n$ tends to a finite endpoint with $\inf_x\rho(t_n,x)^2\geq\mu$,
the local construction gives a common positive existence time about
every $t_n$. Uniqueness then extends the solution across the endpoint,
contradicting maximality. This proves the stated limit at either endpoint.
\end{proof}
\section{Modulation and momentum balance}\label{sec:fitted-modulation}

Write $U=(\eta,\nu)^T$, $Q=Q_\beta=(\hat\eta_\beta,\hat\nu_\beta)^T$,
$c=2\sin\beta$, $P(\beta)=\mathcal P(Q_\beta)$, and
$\sigma=\left(\begin{smallmatrix}0&1\\1&0\end{smallmatrix}\right)$.
A prime on a profile denotes a $y$ derivative. Take $\beta$ in a compact
interval $I$ inside either open gray branch. The profiles and their
derivatives in $y$ and $\beta$ decay exponentially, uniformly on $I$.
Choose $k>0$ small enough that every profile test used below, together
with its required derivatives and polynomial factors, is bounded by
$C\rho^4$, where $\rho(y)=\operatorname{sech}(ky)$.
For $y=x-a(t)$ define
\[
 r(y,t)=U(a(t)+y,t)-Q_{\beta(t)}(y)
       =(\hat\theta,\hat\xi)^T,\qquad
 D_\rho=\int\rho^2(\hat\theta_y^2+\hat\theta^2+\hat\xi^2)\,dy.
\]
With $\chi_R(y)=\operatorname{sech}(y/R)$, let
$D_R=\int\chi_R^2(\hat\theta_y^2+\hat\theta^2+\hat\xi^2)\,dy$.
In particular $D_\rho\leq D_R$ for $R\geq k^{-1}$.

\begin{proposition}[Fitted parameters and their derivatives]\label{prop:mod}
In an energy neighborhood of a translated gray profile the conditions
\begin{equation}\label{eq:fit-conditions}
 \langle r,\sigma Q\rangle=0,\qquad \langle r,Q'\rangle=0
\end{equation}
determine unique nearby parameters $a,\beta$. They are $C^1$ in the datum
and $C^1$ in time on every interval where the solution stays in this
neighborhood. If the datum is initially within $\varepsilon$ of a profile
with parameters $a_0,\beta_0$, then
\[
 |a(0)-a_0|+|\beta(0)-\beta_0|+\|r(0)\|_{\mathcal X}
 \leq C\varepsilon .
\]
Here and below the reference profiles are kept in the interior of the
parameter chart. Put $\alpha=\dot a-c-V(a)$. Then
\begin{align}
 \left|\dot\beta-\frac{V'(a)P(\beta)}{2\cos^2\beta}\right|
 &\leq C(h^3+h\sqrt{D_\rho}+D_\rho),\label{eq:beta-local-fit}\\
 \left|\alpha+
    \frac{2\langle r,\mathcal L_\beta\sigma Q''\rangle}
         {\|Q'\|_2^2}\right|
 &\leq C(h^2+h\sqrt{D_\rho}+D_\rho).\label{eq:alpha-local-fit}
\end{align}
Consequently $|\dot\beta|\leq C(h+D_\rho)$ and
$|\alpha|\leq C(h^2+\sqrt{D_\rho})$. In terms of the dual local energy $D_w$
used in the forced virial estimates, one also has
\begin{equation}\label{eq:alpha-dual-fit}
 |\alpha|\leq C\big(D_w^{1/2}+R^{-2}D_R^{1/2}+D_R+h^2\big).
\end{equation}
The constant in \eqref{eq:alpha-dual-fit} is uniform for $R\geq S\geq S_0(I)$
and $0<\gamma\leq1$, where $S$ is the dual localization radius and
$\gamma$ the regularization parameter.
\end{proposition}

\begin{proof}
The unforced fitted chart with translation and momentum constraints was established for soliton chains by B\'ethuel, Gravejat and Smets \cite[Proposition~6]{BGSChains2014}. We give the argument in the present normalization to obtain the forced derivative estimates. One can put both translations on the smooth profile tests in
\eqref{eq:fit-conditions}; differentiability of translations of the
energy datum itself is unnecessary. In the order $(\beta,a)$ the
Jacobian of the two constraints at $r=0$ has diagonal entries
$-\langle\partial_\beta Q,\sigma Q\rangle=-2P'(\beta)=4\cos^2\beta$
and $\|Q'\|_2^2$. The off-diagonal entries vanish by parity.
The implicit function theorem gives the fitted coordinates and the
initial estimate, uniformly after slightly enlarging the compact
interval $I$. Differentiation of the constraints along an energy
solution can first be performed with the parameters held fixed:
the time derivative of every translated smooth test pairing follows
from the distributional equation and depends continuously on the
datum in $\mathcal X$. The implicit function theorem with time as
an additional parameter thus gives $C^1$ time dependence.

Set $b(y,t)=V(a+y)-V(a)$ and
$\mathcal N=H_0'(Q+r)-H_0'(Q)-H_0''(Q)r$. Since
$J=-2\partial_y\sigma$, $J\mathcal P'(U)=-U_y$ and $(H_0-c\mathcal P)'(Q)=0$, the
equation for $r$ is
\begin{equation}\label{eq:fit-perturbation}
 r_t=J\mathcal L_\beta r+J\mathcal N
   +\alpha(Q'+r_y)-\dot\beta\,\partial_\beta Q+f,\qquad
 f=-\partial_y\big(b(Q+r)\big),
 \quad \mathcal L_\beta=(H_0-c\mathcal P)''(Q).
\end{equation}
The potential appears linearly because $\mathcal P_V(U)=\frac12\int V\eta\nu$
is quadratic. Using $\mathcal L_\beta Q'=0$ and testing
\eqref{eq:fit-perturbation} against the two constraints gives
\begin{align}
 (2P'-\langle r,\sigma\partial_\beta Q\rangle)\dot\beta
     +\langle r,\sigma Q'\rangle\alpha
  &=2\langle\mathcal N,Q'\rangle+\langle f,\sigma Q\rangle ,
    \label{eq:fit-system-one}\\
 \langle r,\partial_\beta Q'\rangle\dot\beta+
   (\|Q'\|_2^2-\langle r,Q''\rangle)\alpha
  &=-2\langle\mathcal L_\beta r+\mathcal N,\sigma Q''\rangle
      -\langle f,Q'\rangle.\label{eq:fit-system-two}
\end{align}
In particular $\langle J\mathcal L_\beta r,\sigma Q\rangle
=2\langle r,\mathcal L_\beta Q'\rangle=0$.

For clarity, the nonlinear pairings in these identities are bounded by
$C D_\rho$ at energy regularity. Indeed the first component of $\mathcal N$
has the form $\partial_y G+F$, and the second is
$-\hat\theta\hat\xi$. Expanding $H_0'$ in divergence form, as in
\cite{HM}, gives
\[
 |G|\leq C(|\hat\theta\hat\theta_y|+\hat\theta^2),\qquad
 |F|\leq C(\hat\theta_y^2+\hat\theta^2+\hat\xi^2).
\]
The denominator $1-\hat\eta_\beta-\hat\theta$ stays bounded away from
zero on the chosen neighborhood. Integrate the divergence against the
derivative of the smooth test to obtain the claimed bound.
All coefficients involving $r$ on the left of
\eqref{eq:fit-system-one}--\eqref{eq:fit-system-two} are $O(\sqrt{D_\rho})$.
The linear Hessian pairing on the right is $O(\sqrt{D_\rho})$ by moving
$\mathcal L_\beta$ onto the localized test.

For either test $g$ we have
$|\langle-\partial_y(br),g\rangle|=|\langle br,g'\rangle|
\leq Ch\sqrt{D_\rho}$, using $|b|\leq h\|W'\|_\infty |y|$.
Write $Q=(Q_1,Q_2)^T$. For the profile part,
\begin{align*}
 \langle-\partial_y(bQ),\sigma Q\rangle
 &=\int b(Q_1Q_2)'=-\int V'(a+y)Q_1Q_2
   =-2V'(a)P(\beta)+O(h^3),\\
 \langle-\partial_y(bQ),Q'\rangle
 &=\int b Q\cdot Q''=O(h^2).
\end{align*}
The first Taylor remainders containing $yQ_1Q_2$ and $yQ\cdot Q''$
integrate to zero by parity. The $O(h^3)$ estimate uses
$\|V'''\|_\infty=O(h^3)$. Invert the matrix in
\eqref{eq:fit-system-one}--\eqref{eq:fit-system-two}, whose diagonal
part at $r=0$ is uniformly invertible. This first bounds both unknowns by
$C(h+\sqrt{D_\rho})$. Substitution into \eqref{eq:fit-system-one} gives
\eqref{eq:beta-local-fit}, since $2P'=-4\cos^2\beta$; substitution into
\eqref{eq:fit-system-two} yields \eqref{eq:alpha-local-fit}.

We verify the cutoff dependence of \eqref{eq:alpha-dual-fit}. Put
$\zeta=1-Q_1$, $\mathsf P=\operatorname{diag}(\sqrt\zeta,1/\sqrt\zeta)$
and $R_\gamma=(1-\gamma\partial_y^2)^{-1}$. The dual variable of
\cite{HM}, also used in the forced virial argument, is
\[
 w=R_\gamma\mathsf P^{-1}\sigma\mathcal L_\beta(\chi_R r),
 \qquad D_w=\int\chi_S^2
    (|\partial_y w_1|^2+|w_1|^2+|w_2|^2)\,dy .
\]
In the principal pairing of \eqref{eq:alpha-local-fit}, insert the
broad cutoff before applying the Hessian. Distributional self-adjointness
gives
\begin{align*}
 \langle r,\mathcal L_\beta\sigma Q''\rangle
 &=\left\langle w,(1-\gamma\partial_y^2)
                  \left(\frac{\mathsf P Q''}{\chi_R}\right)\right\rangle
   -\left\langle[\mathcal L_\beta,\chi_R]r,
                    \frac{\sigma Q''}{\chi_R}\right\rangle .
\end{align*}
The divided smooth test and its first two derivatives have bounded
$L^2(\chi_S^{-2}\,dy)$ norm, uniformly for $R\geq S\geq S_0(I)$
by exponential profile decay. The first term is at most $C D_w^{1/2}$.
With $m=\zeta^{-1}$ the only nonzero component of the commutator is
$-\frac14(2m\chi_R'\hat\theta_y+(m\chi_R')'\hat\theta)$.
The bounds $|\chi_R'|/\chi_R\leq |y|/R^2$ and
$|\chi_R''|/\chi_R\leq C/R^2$ give the second term a bound
$CR^{-2}D_R^{1/2}$ by moving the remaining inverse cutoff onto the
exponentially decaying test. Applying $h\sqrt{D_\rho}\leq(h^2+D_\rho)/2$ in
\eqref{eq:alpha-local-fit} completes the proof.
\end{proof}

\begin{proposition}[Exact momentum balance and its force defect]\label{prop:force}
On any interval of the preceding decomposition let
$\Pi(t)=\mathcal P(U(t))$. Then
\begin{align}
 \Pi(t)&=P(\beta(t))+\tfrac12\int\hat\theta\hat\xi\,dy,
 &\dot\Pi(t)&=-\tfrac12\int V'(x)\eta(x,t)\nu(x,t)\,dx ,
 \label{eq:fit-total-momentum}\\
 |\dot\Pi(t)+V'(a(t))\Pi(t)|
 &\leq C\big(h^3+h^2\sqrt{D_\rho}+h\|r\|_{\mathcal X}^2\big).
 \label{eq:fit-force-defect}
\end{align}
In particular, the leading slow momentum equation is
$dP(B)/ds=-W'(A)P(B)$, with $s=ht$, $A=ha$.
\end{proposition}

\begin{proof}
The first identity follows by expanding the quadratic $\mathcal P$ and using
\eqref{eq:fit-conditions}. The second is the energy-space balance \eqref{eq:momentum}.

In moving coordinates insert the factor $V'(a+y)-V'(a)$ in this last
integral. In the pure profile product subtract also $yV''(a)$ by parity;
the resulting bound is $Ch^3\int y^2|Q_1Q_2|$. The two terms linear
in $r$ cost $Ch^2\sqrt{D_\rho}$, since the coefficient difference is at most
$Ch^2|y|$ and the other factor is a profile. In the remaining term use
$|V'(a+y)-V'(a)|\leq 2h\|W'\|_\infty$ and
$\|\hat\theta\hat\xi\|_1\leq\|r\|_{\mathcal X}^2/2$.
This gives \eqref{eq:fit-force-defect} for radiation located anywhere
on the line.
\end{proof}

For a comparison trajectory initialized as $(A(0),B(0))=(ha(0),\beta(0))$,
\eqref{eq:fit-total-momentum} gives
$|\Pi(0)-P(B(0))|\leq C\|r(0)\|_{\mathcal X}^2$.
Starting the effective flow at an arbitrary reference profile within
$O(\sqrt h)$ of the data can produce an initial speed-parameter discrepancy
of order $\sqrt h$. This distinction already occurs for neighboring exact
gray solitons with $V=0$.
\section{The modulated action under transport}\label{sec:action}

Write $Z_\beta=H_0-2\sin\beta\,\mathcal P$ and
\[
 Y(t)=Z_{\beta(t)}(U(t))-Z_{\beta(t)}(Q_{\beta(t)}).
\]
The two terms can be evaluated in the frame $y=x-a(t)$ because
$Z_\beta$ is translation invariant. Throughout this section $\beta$
lies in the fixed compact gray interval $I$. We use $D_R$ for the
localized squared energy of $r=(\hat\theta,\hat\xi)^T$ with the weight
$\operatorname{sech}^2(y/R)$, taking $R$ sufficiently large in terms of $I$.

\begin{proposition}[Modulated action]\label{prop:action}
On the fitted neighborhood of Proposition~\ref{prop:mod}, for
$\|r\|_{\mathcal X}$ sufficiently small there are $c_I,C_I>0$ such that
\begin{equation}\label{eq:action-equiv}
 c_I\|r\|_{\mathcal X}^2\le Y\le C_I\|r\|_{\mathcal X}^2 .
\end{equation}
The action obeys
\begin{equation}\label{eq:action-rate}
 |Y'|\le C\left[hD_R^{1/2}
       +h\|r\|_{\mathcal X}^2+|\dot\beta|\,\|r\|_{\mathcal X}^2\right].
\end{equation}
The constants in the derivative bound may also depend on $W$. These
statements require only $H^1\times L^2$ control of the perturbation.
\end{proposition}

\begin{proof}
The critical point identity $Z'_\beta(Q_\beta)=0$ gives
\[
 Y=\tfrac12\langle\mathcal L_\beta r,r\rangle+\mathcal R,\qquad
 \mathcal L_\beta=Z''_\beta(Q_\beta).
\]
To see the size of $\mathcal R$, put $\zeta=1-\hat\eta_\beta$. Direct
subtraction of the quadratic action gives
\[
 \mathcal R=\int\left[
 \frac{\hat\theta}{8\zeta(\zeta-\hat\theta)}
 \left(\hat\theta_y+
       \frac{\hat\eta_\beta'}{\zeta}\hat\theta\right)^2
 -\tfrac12\hat\theta\hat\xi^2\right]dy .
\]
Both $\zeta$ and $\zeta-\hat\theta$ have a uniform positive lower
bound, since $\inf_y\zeta=\sin^2\beta$ and $H^1(\mathbb R)$ embeds
into $L^\infty$. Thus $|\mathcal R|\le C_I\|r\|_{\mathcal X}^3$.
B\'ethuel, Gravejat and Smets proved constrained coercivity for the single gray soliton \cite[Proposition~4 and equations~(21)--(23)]{BGSChains2014}. In our normalization the estimate of \cite[Section~2]{HM} applies to
the fitted conditions $\langle r,\sigma Q_\beta\rangle
=\langle r,\partial_yQ_\beta\rangle=0$ and gives
$\langle\mathcal L_\beta r,r\rangle\ge c\|r\|_{\mathcal X}^2$,
uniformly on $I$. The form is bounded above on $\mathcal X$, which
proves \eqref{eq:action-equiv}.

We compute $Y'$ first on smooth solutions. The unforced vector
field preserves both $H_0$ and $\mathcal P$, and hence $Z_\beta$ at
fixed $\beta$. In moving coordinates the additional term
$\dot a\,\partial_y U$ pairs to zero with $Z'_\beta(U)$ by translation
invariance. The derivative of the profile in
$Z_{\beta(t)}(Q_{\beta(t)})$ pairs to zero with $Z'_\beta(Q_\beta)$.
The remaining contributions are the transport by
$b(y)=V(a+y)-V(a)$ and the explicit angle derivative. In particular,
the latter is exactly
\begin{equation}\label{eq:action-angle}
 -2\cos\beta\,\dot\beta\,[\mathcal P(U)-P(\beta)]
 =-\cos\beta\,\dot\beta\int\hat\theta\hat\xi\,dy
\end{equation}
by Proposition~\ref{prop:mod}.

For the transport calculation, write $(e,d,v)=(\eta,\eta_y,\nu)$
in the moving coordinates and let $z$ denote the density of $Z_\beta$,
\[
 z(e,d,v)=\frac{d^2}{8(1-e)}+\frac12(1-e)v^2
             +\frac12 e^2-\sin\beta\,ev .
\]
Under $(e,v)_t=-\partial_y[b(e,v)]$, the chain rule gives
$z_t=-b\partial_yz-b'(ez_e+2dz_d+vz_v)-b''ez_d$.
Integrating in $y$ gives the transport contribution
$\int(b'S_\beta+b''T)\,dy$, with
\begin{align*}
 S_\beta(e,d,v)
   &=z-ez_e-2dz_d-vz_v\\
   &=-\frac{3-2e}{8(1-e)^2}d^2
      +\left(e-\tfrac12\right)v^2-\tfrac12 e^2+\sin\beta\,ev,\\
 T(e,d,v)&=-ez_d=-\frac{e\,d}{4(1-e)} .
\end{align*}
For a pair $F$ write $S_\beta[F]=S_\beta(F_1,\partial_y F_1,F_2)$
and similarly for $T[F]$. At $Q_\beta$ the integral of the two
transport terms is zero since
$\langle Z'_\beta(Q_\beta),-\partial_y(bQ_\beta)\rangle=0$.
Thus the exact stress identity, with $G=Q_\beta+r$, is
\begin{align}
 Y'={}&\int b'\big(S_\beta[G]-S_\beta[Q_\beta]\big)\,dy
       +\int b''\big(T[G]-T[Q_\beta]\big)\,dy\notag\\
      &-2\cos\beta\,\dot\beta\,[\mathcal P(U)-P(\beta)] .
 \label{eq:action-stress}
\end{align}
Only the spatial derivatives $b'=V_x(a+y)$ and $b''=V_{xx}(a+y)$
occur.

To bound the differences, set $q=(\hat\eta_\beta,\hat\eta'_\beta,
\hat\nu_\beta)$ and $v_r=(\hat\theta,\hat\theta_y,\hat\xi)$.
For $F=S_\beta,T$, the explicit formulas yield
\[
 |\nabla F(q)|\le C_I e^{-\gamma_I|y|},\qquad
 |F(q+v_r)-F(q)-\nabla F(q)\cdot v_r|
     \le C_I(\hat\theta^2+\hat\theta_y^2+\hat\xi^2).
\]
For example, for a term $m(e)d^2$ the remainder equals
\[
 m(e_0+\hat\theta)\hat\theta_y^2
 +2[m(e_0+\hat\theta)-m(e_0)]d_0\hat\theta_y
 +[m(e_0+\hat\theta)-m(e_0)-m'(e_0)\hat\theta]d_0^2,
\]
where $(e_0,d_0)=(\hat\eta_\beta,\hat\eta'_\beta)$.
The coefficient and its first two derivatives are bounded on the
range in question. The $v^2$ term has the same expansion with
$\hat\xi$ in place of $\hat\theta_y$. For $T=m(e)d$ one uses
$|m(e_0+\hat\theta)-m(e_0)|\le C_I|\hat\theta|$ in the mixed
term. The bilinear and quadratic terms in $e,v$ satisfy the
displayed bound directly. By weighted Cauchy--Schwarz the integrated
linear terms cost $C_I D_R^{1/2}$ for $R>1/\gamma_I$. Therefore
the two stress integrals in \eqref{eq:action-stress} have absolute sum
at most
\[
 C_I(\|V_x\|_\infty+\|V_{xx}\|_\infty)
            (D_R^{1/2}+\|r\|_{\mathcal X}^2).
\]
Now $\|V_x\|_\infty=h\|W'\|_\infty$ and
$\|V_{xx}\|_\infty=h^2\|W''\|_\infty$. Combining with
\eqref{eq:action-angle} proves \eqref{eq:action-rate} for $h\le1$.

Finally, the integrations by parts for smooth solutions can be made
with expanding spatial cutoffs, since the densities and stresses
are integrable and $b,b',b''$ are bounded. On a finite interval in
the fitted neighborhood, approximate the solution using
Proposition~\ref{prop:flow}. Strong convergence in
$C_t(H^1\times L^2)$ and the uniform lower bound on $1-\eta$ imply
convergence in $C_tL^1_y$ of the two stresses. Indeed the density
and its coefficients converge uniformly, while products of two
strongly convergent $L^2$ factors converge in $L^1$.
The moving centers and angles converge with their fitted equations
by Proposition~\ref{prop:mod}. Translations act continuously in
$\mathcal X$, and $V_x,V_{xx}$ depend continuously in the supremum
norm on the shift $a$. Passing to the integrated form of
\eqref{eq:action-stress} proves the identity, hence the derivative
bound, also at energy regularity.
\end{proof}
\begin{corollary}[Energy of the fitted parameters]\label{cor:fitted-energy}
Under the hypotheses of Proposition~\ref{prop:action}, let
$e_{\rm fit}(t)=\HH(ha(t),\beta(t))$. For $0<h\le1$,
\begin{align}
 |H(U(t))-e_{\rm fit}(t)|
 &\le C_{I,W}\big(h^2+\|r(t)\|_{\X}^2\big),\label{eq:fitted-energy-level}\\
 |e_{\rm fit}(t)-e_{\rm fit}(0)|
 &\le C_{I,W}\big(h^2+\|r(0)\|_{\X}^2+\|r(t)\|_{\X}^2\big).
 \label{eq:fitted-energy-drift}
\end{align}
\end{corollary}
\begin{proof}
Write $Q=Q_{\beta(t)}$, $r=(r_1,r_2)^T$, $c=2\sin\beta(t)$, and
$b(y)=V(a+y)-V(a)$. The definitions give
\[
 H(U)-e_{\rm fit}
 =Y+(c+V(a))[\PP(U)-P(\beta)]
   +\tfrac12\int b\,(Q_1Q_2+Q_1r_2+Q_2r_1+r_1r_2)\,dy .
\]
Propositions~\ref{prop:force} and~\ref{prop:action} bound the first
two terms by $C_{I,W}\|r\|_{\X}^2$. In the pure profile term,
subtracting $yV'(a)$ from $b$ leaves the integral unchanged by parity,
and the Taylor remainder costs $C_{I,W}h^2$.
The mixed terms are bounded by $C_{I,W}h\|r\|_{\X}$ since
$|b|\le h\|W'\|_\infty|y|$ and $yQ\in L^2$ uniformly on $I$.
The last term costs $C_W\|r\|_{\X}^2$ since $|b|\le2\|W\|_\infty$.
Young's inequality proves \eqref{eq:fitted-energy-level}.
Conservation of $H$ gives \eqref{eq:fitted-energy-drift}.
\end{proof}

\subsection{The second moment of the core momentum}
With $M_2$ as in \eqref{eq:finite-width-hamiltonian}, one has
\[
 M_2(\beta)=
 \frac{(\pi^2/4)P(\beta)-(\pi/2-\beta)^3}{3\cos^2\beta}.
\]
To evaluate it, write $d=\cos\beta$ and use
$2p_\beta=\hat\nu_\beta-\sin\beta\,\hat\eta_\beta$.
For $y>0$ the geometric series gives
$\hat\nu_\beta(y)=2d\sum_{n\ge1}e^{-2ndy}\sin(2n\beta)$.
Integrating termwise against $y^2$ and using
$\int_0^\infty x^2\sech^2 x\,dx=\pi^2/12$ yields
\[
 M_2(\beta)=\frac{1}{2d^2}\sum_{n\ge1}\frac{\sin(2n\beta)}{n^3}
              -\frac{\pi^2\sin\beta}{12d}.
\]
The sine sum at $z=2\beta$ equals
$\pi^2 z/6-\pi z^2/4+z^3/12$ for $0\le z\le2\pi$, by integrating
the Fourier series for $(z-\pi)^2$. This gives the displayed formula.
For $U(0)=Q_\beta(\cdot-a)$ the momentum balance \eqref{eq:momentum}
gives the initial force
\[
 \left.\frac{d}{dt}\PP(U(t))\right|_{t=0}
 =-\partial_a H(Q_\beta(\cdot-a))
 =-hW'(ha)P(\beta)-\tfrac12h^3 W'''(ha)M_2(\beta)+O_{I,W}(h^5).
\]
The second moment gives the initial force at this order. Higher-order tracking of the fitted center would also require the radiation and the particular translation constraint in \eqref{intro:orth}.

\section{Local virial estimates in the momentum field}
\label{vir:section}

B\'ethuel, Gravejat and Smets applied the Hessian transform and a linear-weight virial to a limiting Gross--Pitaevskii orbit \cite[equation~(34) and Proposition~7]{BGS_Asymp_Dark_2015}. Here we use the fixed-time spatial estimates of \cite{HM} for the evolving fitted perturbation. The forced terms below require control of radiation transported through the region where the field varies.

Write $r=(\hat\theta,\hat\xi)^T$ in the frame $y=x-a(t)$,
$\zeta_\beta=1-\hat\eta_\beta$, $b(y)=V(a+y)-V(a)$, and
$\alpha=\dot a-2\sin\beta-V(a)$. Throughout this section
$\beta\in I$, where $I$ is a compact interval in a gray branch.
For $L\geq1$ put
\begin{gather*}
 \chi_L(y)=\operatorname{sech}(y/L),\qquad
 \phi_L(y)=L\tanh(y/L),\qquad m_L=\zeta_\beta\phi_L,\\
 D_L(f)=\int\chi_L^2\bigl(|f_1'|^2+|f_1|^2+|f_2|^2\bigr).
\end{gather*}
We abbreviate $D_R=D_R(r)$, with $R$ a fixed broad radius. Fix a
second, narrower radius $S$ and a regularization parameter $\gamma>0$.
Using the Hessian $\mathcal L_\beta$ of $H_0-2\sin\beta\,\mathcal P$
at $Q_\beta$, define
\begin{gather*}
 \mathsf P_\beta=\operatorname{diag}(\sqrt{\zeta_\beta},1/\sqrt{\zeta_\beta}),
 \qquad \mathcal T_\beta=\mathsf P_\beta^{-1}\sigma\mathcal L_\beta,
 \qquad R_\gamma=(1-\gamma\partial_y^2)^{-1},\\
 w=R_\gamma\mathcal T_\beta(\chi_Rr),\qquad D_w=D_S(w).
\end{gather*}
The first row of $\mathcal T_\beta$ consists of multiplication operators.
The second row contains a second-order operator on the density and a
multiplication operator on the velocity. These coefficients have bounded
derivatives uniformly on $I$. Consequently
\begin{equation}\label{vir:wb}
 \|w_1\|_{H^2}+\|w_2\|_{H^1}\leq C_\gamma D_R^{1/2}.
\end{equation}
The primal and dual functionals are
\[
 \mathcal I_R=\int m_R\hat\theta\hat\xi,\qquad
 \mathcal V_{S,\gamma}=\int m_S w_1w_2 .
\]

\begin{proposition}[Forced local virial estimate]\label{prop:vir}
There are fixed $S,\gamma,R,M,c,C>0$, a radius $\delta_0>0$ and $h_0>0$
such that, whenever $\|r\|_{\mathcal X}\leq \delta_0$, $0<h\leq h_0$ and
$\beta\in I$, the functional $\mathcal F=\mathcal I_R+M\mathcal V_{S,\gamma}$
satisfies
\begin{equation}\label{vir:combined}
 |\mathcal F|\leq C\|r\|_{\mathcal X}^2,\qquad
 \mathcal F'\geq c(D_R+D_w)-C\bigl(h^2+h\|r\|_{\mathcal X}^2\bigr).
\end{equation}
The differential inequality holds in integrated form for energy solutions.
The constants are independent of the length of the time interval.
\end{proposition}

\subsection{The primal estimate}
Start with a smooth solution and write $\delta=\|r\|_{\mathcal X}$.
The static density recovery estimate \cite[Lemma~4.3]{HM} gives
\begin{equation}\label{vir:density}
 \|\rho\hat\theta\|_2^2\leq C D_w+C R^{-4}D_R .
\end{equation}
The forced parameter estimates in Proposition~\ref{prop:mod} imply
\begin{align}
 |\dot\beta|&\leq C(h+D_R),\qquad
 |\alpha|\leq C(h^2+D_R^{1/2}),\label{vir:params}\\
 |\alpha|&\leq C\bigl(h^2+D_w^{1/2}+R^{-2}D_R^{1/2}+D_R\bigr).
 \label{vir:refined}
\end{align}
All three bounds are uniform for $R\geq S$ sufficiently large and
$0<\gamma\leq1$. The $D_w^{1/2}$ bound on the linear translation
functional will allow us to absorb it by the dual virial.

The linear and nonlinear parts $J\mathcal L_\beta r+J\mathcal N$
of \eqref{eq:fit-perturbation} contribute
$2\langle\mathcal L_\beta r+\mathcal N,\partial_y(m_Rr)\rangle$
to $\mathcal I_R'$. The primal linear form and nonlinear
estimate of \cite[Lemmas~5.1--5.2]{HM} give the lower bound
\[
 c_0D_R-C\|\rho\hat\theta\|_2^2-C\delta D_R.
\]
These are spatial form inequalities for the given $r,\beta$.
The transport terms $\alpha r_y-\partial_y(br)$ have exactly the
contribution
\[
 -\alpha\int m_R'\hat\theta\hat\xi
 +\int (b m_R'-b' m_R)\hat\theta\hat\xi .
\]
Since $|m_R'|\leq C\chi_R^2$, the first term is
$O((h+\delta)D_R)$. For the second term both coefficients are bounded
by $C_R h$ using $|b|\leq Ch|y|$, the decay of $m_R'$, and
$\|m_R\|_\infty\leq R$. Its global bound $C_R h \delta^2$ permits $b$
to have order-one size far from the soliton.

The localized profile force $-(bQ_\beta)'$ pairs with $m_R\sigma r$
at a cost $C h D_R^{1/2}$. Indeed $|m_R|\leq |y|$ and every term
in this source has an exponentially decaying coefficient bounded
by $C h(1+|y|)e^{-c_I|y|}$. The $\partial_\beta Q_\beta$ term
costs $C(h+D_R)D_R^{1/2}$, and differentiating $m_R$ in $\beta$
costs $C(h+D_R)D_R$. For the translation of the profile, the bound is
$C|\alpha|D_R^{1/2}$. By \eqref{vir:refined}, for any $\varepsilon>0$
this is at most
\[
 2\varepsilon D_R+C_\varepsilon(D_w+h^2)
       +C(R^{-2}+\delta)D_R.
\]
Use \eqref{vir:density} in the linear form, choose $\varepsilon$ small,
then take $R$ sufficiently large and $\delta,h$ sufficiently small.
Together the bounds give
\begin{equation}\label{vir:primal}
 \mathcal I_R'\geq c_1D_R-C_1D_w-C_R(h^2+h \delta^2).
\end{equation}
Here $c_1,C_1$ are independent of $R$ large and $\gamma\leq1$.

\subsection{The regularized dual estimate}
Write $\mathcal T=\mathcal T_\beta$. The frozen dual generator is
$\mathcal G_\beta=-2\mathsf P_\beta^{-1}\sigma\mathcal L_\beta
\partial_y\mathsf P_\beta$, with the derivative acting by
composition. Thus $\mathcal T J\mathcal L_\beta=\mathcal G_\beta\mathcal T$.
Differentiation gives
\begin{align*}
 w_t={}&\mathcal G_\beta w+[R_\gamma,\mathcal G_\beta]R_\gamma^{-1}w
       +R_\gamma\mathcal T\bigl([\chi_R,J\mathcal L_\beta]r
                                     +\chi_R J\mathcal N\bigr)\\
 &+\alpha w_y-\alpha R_\gamma
       \bigl((\partial_y\mathcal T)(\chi_Rr)+\mathcal T(\chi_R'r)\bigr)
       +\dot\beta R_\gamma(\partial_\beta\mathcal T)(\chi_R r)\\
 &+\alpha R_\gamma\mathcal T(\chi_R Q_\beta')
       -\dot\beta R_\gamma\mathcal T(\chi_R\partial_\beta Q_\beta)
       -R_\gamma\mathcal T\bigl(\chi_R(bQ_\beta)'\bigr)
       -R_\gamma\mathcal T\bigl(\chi_R(br)'\bigr).
\end{align*}
Here $\partial_y\mathcal T$ differentiates coefficients only.
We record the bounds for these terms in their virial pairing against
$m_S\sigma w$. For any source $f$ one has
\begin{equation}\label{vir:pair}
 |\langle f,m_S\sigma w\rangle|
 \leq C_S(\|f_1\|_2+\|f_2\|_{H^{-1}})\|w\|_{\mathcal X}
 \leq C_{S,\gamma}(\|f_1\|_2+\|f_2\|_{H^{-1}})D_R^{1/2}.
\end{equation}

The frozen generator contributes the quadratic form of
$\partial_y\mathcal L_\beta m_S-m_S\mathcal L_\beta\partial_y$
on $\mathsf P_\beta w$, since the diagonal entries of $\mathsf P_\beta$
have product one. The truncated coercivity and regularization
commutator estimates \cite[Lemmas~6.1--6.2]{HM} bound the first two
terms from below by
$\kappa_0D_w-CR^{-4}D_R-C_S\sqrt\gamma D_w$.
Here the constrained coercivity is applicable because
\[
 \big\langle w,R_\gamma^{-1}\mathsf P_\beta\sigma\partial_\beta Q_\beta\big\rangle
   =\cos\beta\,\langle\chi_Rr,\sigma Q_\beta\rangle
   =O(R^{-2}D_R^{1/2}).
\]
The last step follows from orthogonality and
$\|(1-\chi_R^{-1})\sigma Q_\beta\|_2\leq C R^{-2}$.
This is the fixed-time transfer in \cite[Lemma~4.2]{HM}.
The two sources remaining on the first line of $w_t$ are bounded
in $L^2\times H^{-1}$ by $C_\gamma (R^{-1}+\delta)D_R^{1/2}$, by
\cite[Lemma~6.3]{HM}. The weighted momentum is quadratic, so $\mathcal N$ contains only terms from $H_0$. Formula \eqref{vir:pair} gives their
contribution $O(C_{S,\gamma}(R^{-1}+\delta)D_R)$.

The advection term pairs to $-\alpha\int m_S'w_1w_2$, of size
$C(h+\delta)D_w$. The remaining terms on the second line and the
time derivative of $m_S$ have total size at most
$C_{S,\gamma}(|\alpha|+|\dot\beta|)D_R$. Indeed the row orders of
$\partial_y\mathcal T$ and $\partial_\beta\mathcal T$ are bounded
by those of $\mathcal T$, $|\partial_y^j\chi_R|\leq C_j R^{-j}\chi_R$,
and \eqref{vir:wb} applies with these coefficient changes. For
the translated profile use $\mathcal T Q_\beta'=0$: the smooth
vector $\mathcal T(\chi_R Q_\beta')$ equals
$[\mathcal T,\chi_R]Q_\beta'$ and has $L^2\times L^2$ norm
$O(R^{-2})$. Its pairing costs
$C_{S,\gamma}R^{-2}(h+D_R^{1/2})D_R^{1/2}$ by \eqref{vir:params}.
The $\partial_\beta Q_\beta$ and $(bQ_\beta)'$ sources have respective
$L^2\times H^{-1}$ norms after transformation at most
$C_\gamma(h+D_R)$ and $C_\gamma h$, by profile decay and the
bounds for the derivatives of $b$. Their pairings therefore cost
$C_{S,\gamma}(hD_R^{1/2}+\delta D_R)$.

It remains to control transport of the radiation. For fixed $R$,
\[
 \|\chi_R b\|_{W^{1,\infty}}+\|\chi_R'b\|_\infty
       +\|\chi_R b'\|_\infty\leq C_R h.
\]
Let $g=\chi_R(br)'$. Its first component is
$(\chi_R b)\hat\theta'+\chi_Rb'\hat\theta$; the second equals
$\partial_y((\chi_Rb)\hat\xi)-\chi_R'b\hat\xi$. Hence
$\|g_1\|_2+\|g_2\|_{H^{-1}}\leq C_R h \delta$.
The row structure of $\mathcal T$ gives a bounded map
$L^2\times H^{-1}\to H^{-1}\times H^{-2}$ and the resolvent
gains two derivatives. In particular,
\begin{equation}\label{vir:radiation}
 \|(R_\gamma\mathcal T g)_1\|_{H^1}
       +\|(R_\gamma\mathcal T g)_2\|_2\leq C_{R,\gamma}h \delta.
\end{equation}
Thus the virial pairing costs $C_{R,S,\gamma}h \delta D_R^{1/2}$.
This argument puts the sole derivative $\partial_y\hat\xi$ inside a distributional derivative of an $L^2$
function.

Choose $S$ large and $\gamma$ small as in the frozen inequalities.
For every $\varepsilon>0$, Young's inequality on the terms containing
$hD_R^{1/2}$ and $h\delta D_R^{1/2}$ then gives
\begin{equation}\label{vir:dual}
 \mathcal V_{S,\gamma}'\geq
 \kappa_1D_w-\bigl[\varepsilon+C_{S,\gamma}(R^{-1}+\delta+h)\bigr]D_R
       -C_{R,S,\gamma,\varepsilon}h^2(1+\delta^2),
\end{equation}
where $\delta,h$ are small enough to absorb the advection loss in $D_w$.
In particular the constants depending on $R$ in \eqref{vir:radiation}
multiply a term of order $h^2$ in \eqref{vir:dual}.

\begin{proof}[Proof of Proposition~\ref{prop:vir}]
Fix $M$ with $M\kappa_1>C_1+1$ in \eqref{vir:primal}--\eqref{vir:dual},
then $\varepsilon$ with $M\varepsilon<c_1/4$. Increase $R$ and
decrease $\delta_0,h_0$ to absorb the remaining $D_R$ errors. The
size bound follows from boundedness of the two weights and
\eqref{vir:wb}. This proves \eqref{vir:combined} for smooth solutions.
On any compact time interval on which the solution lies in the
nonvanishing modulation neighborhood, use smooth data approximation
in $C_t\mathcal X$ and the continuity of the modulation parameters.
The operators defining $w$ are continuous in $\beta$ from
$\mathcal X$ to $H^2\times H^1$, uniformly on $I$ at fixed cutoffs.
Both endpoint functionals, both local energies, and the error
integrals consequently converge. Passing in the integrated smooth
inequality proves the assertion at energy regularity.
\end{proof}
\section{Orbital control, trajectory comparison, and first exit}
\label{sec:closure}

All estimates until the continuation argument are on a finite interval
on which the fitted parameter belongs to $I$ and the remainder has
sufficiently small $\mathcal X$ norm. Write
$d_{\rm in}=\|r(0)\|_{\mathcal X}$ and $D=D_R+D_w$. The localization radii
and the regularization parameter in Proposition~\ref{prop:vir} are
fixed independently of $h$ and of the interval length.

\begin{proposition}[A priori estimate on a finite interval]
\label{prop:combined}
There are constants $C,h_0>0$ and a fixed small tube radius such that,
as long as the above conditions hold and $0<h<h_0$,
\begin{equation}\label{eq:combined-finite}
 \sup_{0\le u\le t}\|r(u)\|_{\mathcal X}^2+\int_0^t D(u)\,du
 \le C(d_{\rm in}^2+h^2t)e^{Cht}.
\end{equation}
\end{proposition}

\begin{proof}
Let $Y=Z_\beta(U)-Z_\beta(Q_\beta)$ and let $\mathcal F$ be the combined
virial of Proposition~\ref{prop:vir}. By constrained coercivity,
$c_0\|r\|_{\mathcal X}^2\le Y\le C_0\|r\|_{\mathcal X}^2$.
The virial estimate and its endpoint bound give
\[
 |\mathcal F|\le C_F Y,\qquad \mathcal F'\ge c_1D-C(hY+h^2).
\]
For each fixed $\tau>0$, Proposition~\ref{prop:action} gives
\[
 Y'\le \tau D+C_\tau(hY+h^2)
\]
on decreasing the tube radius if necessary. Indeed the localized
linear error $Ch\sqrt{D_R}$ is at most
$\tfrac12\tau D_R+C_\tau h^2$. A term $CYD_R$ coming from the
parameter derivative, using $|\dot\beta|\le C(h+D_R)$ from Proposition~\ref{prop:mod}, is at most $\tfrac12\tau D_R$ in this tube.
Choose $\lambda>0$ with $\lambda C_F\le 1/2$, followed by
$\tau=\lambda c_1/2$. Then $Y/2\le Y-\lambda\mathcal F\le 3Y/2$, and
\[
 (Y-\lambda\mathcal F)' +\tfrac12\lambda c_1 D
 \le C h(Y-\lambda\mathcal F)+C h^2 .
\]
Multiplying by $e^{-Cht}$ and integrating controls both
$e^{-Cht}(Y-\lambda\mathcal F)(t)$ and
$\int_0^t e^{-Chu}D(u)\,du$ by $C(d_{\rm in}^2+h^2t)$.
This proves \eqref{eq:combined-finite}. The same argument applies to
the integrated differential inequalities obtained by smooth
approximation on each finite existence interval.
\end{proof}

\subsection{Comparison in position and total momentum}
Let $X(t)=ha(t)$, $p(t)=\mathcal P(U(t))$ and $v(\beta)=2\sin\beta$.
In the estimates below $s=ht$. The exact orthogonality identity gives
\begin{equation}\label{eq:compare-fitted-momentum}
 |p(t)-P(\beta(t))|
 =\left|\tfrac12\int\hat\theta\hat\xi\right|
 \le C\|r(t)\|_{\mathcal X}^2 .
\end{equation}
In particular, even when the reduced orbit starts at the fitted
parameters, its initial momentum $P(B(0))$ can differ from $p(0)$
by $O(d_{\rm in}^2)$.
By Propositions~\ref{prop:force} and~\ref{prop:mod},
\begin{align}
 \dot X &= h\big(v(\beta)+W(X)\big)+e_X,
 & |e_X|&\le C h\big(h+\sqrt{D_R}\big),
 \label{eq:center-defect-integral}\\
 \dot p &=-h W'(X)p+e_p,
 & |e_p|&\le C\big(h^3+h^2\sqrt{D_R}
                     +h\|r\|_{\mathcal X}^2\big).
 \label{eq:momentum-defect-integral}
\end{align}
Set $m=d_{\rm in}^2+h$. From \eqref{eq:combined-finite} we have
\begin{equation}\label{eq:comparison-budget}
 \sup_{u\le t}\|r(u)\|_{\mathcal X}^2+\int_0^t D
 \le C m(1+s)e^{Cs}.
\end{equation}
Cauchy--Schwarz in physical time therefore yields
\[
 h\int_0^t\sqrt{D_R(u)}\,du
 \le C\sqrt{h m}\,\sqrt{s(1+s)}e^{Cs/2}
 \le C m e^{C's}.
\]
The latter inequality uses $h\le m$. Also
$h\int_0^t\|r\|_{\mathcal X}^2\le C m s(1+s)e^{Cs}$.
Together with $h^2t=hs$ these inequalities give, on enlarging the
exponent constant,
\begin{equation}\label{eq:comparison-integral}
 \int_0^t\big(|e_X|+|e_p|+h\|r\|_{\mathcal X}^2\big)\,du
 \le C m e^{Cs}.
\end{equation}
The single time integral of the localized linear center error is
where the spacetime part of \eqref{eq:combined-finite} is needed.

Let $(A,B)$ solve the reduced system with
$(A(0),B(0))=(X(0),\beta(0))$. While $B(ht)$ belongs to $K$ and
$\beta(t)$ to $I$, the identity $P'=-2\cos^2\beta$ and
\eqref{eq:compare-fitted-momentum} imply
\[
 |\beta(t)-B(ht)|
 \le C_I\big(|p(t)-P(B(ht))|+\|r(t)\|_{\mathcal X}^2\big).
\]
Define $L(t)=|X(t)-A(ht)|+|p(t)-P(B(ht))|$. Subtracting the reduced
equations in physical time from
\eqref{eq:center-defect-integral}--\eqref{eq:momentum-defect-integral},
and using boundedness of $W'$, $W''$, and $P(B)$, gives
\[
 L(t)\le C d_{\rm in}^2+C h\int_0^t L(u)\,du
 +C\int_0^t\big(|e_X|+|e_p|+h\|r\|_{\mathcal X}^2\big)\,du .
\]
Since the last integral has an increasing upper bound
\eqref{eq:comparison-integral}, Gronwall's inequality and
\eqref{eq:comparison-budget} yield
\begin{equation}\label{eq:trajectory-closed}
 |ha(t)-A(ht)|+|\beta(t)-B(ht)|
 \le C(d_{\rm in}^2+h)e^{Cht}.
\end{equation}
This estimate uses only the initial distance to the \emph{fitted}
profile. Starting the reduced system at an arbitrary reference profile
would introduce its initial parameter discrepancy into the comparison.

\subsection{Direct comparison of the fitted parameters}
Let $F$ be the vector field for position and angle in \eqref{ode}, so
$(A',B')=F(A,B)$. On the same fitted interval, Proposition~\ref{prop:mod} gives
\[
 \left|\frac{d}{dt}(X,\beta)-hF(X,\beta)\right|
 \le C_{I,W}\big(h^3+h\sqrt{D_R}+D_R\big).
\]
Indeed the center defect is $h\alpha$, and the angle defect is
estimated by \eqref{eq:beta-local-fit}. The map $F$ is uniformly
Lipschitz on $\R\times I$. The fitted and reduced parameters agree at
$t=0$, so Gronwall's inequality gives, with $J(t)=\int_0^t D_R(u)\,du$,
\[
 |X(t)-A(ht)|+|\beta(t)-B(ht)|
 \le C e^{Cht}\big(h^3t+h\sqrt{tJ(t)}+J(t)\big).
\]
Using $2h\sqrt{tJ(t)}\le h^2t+J(t)$ and
\eqref{eq:combined-finite}, and increasing $C$, we obtain
\[
 |X(t)-A(ht)|+|\beta(t)-B(ht)|
 \le C(d_{\rm in}^2+h^2t)e^{Cht}.
\]
In particular exact soliton data satisfy
\begin{equation}\label{eq:exact-tracking}
 |ha(t)-A(ht)|+|\beta(t)-B(ht)|\le C h^2t\,e^{Cht}.
\end{equation}
The constants have the same uniformity in the initial center and the
length of the interval as in Proposition~\ref{prop:combined}.

\begin{proof}[Conclusion of the proof of Theorem~\ref{thm:main}]
Suppose $d_{\rm in}\le C_{\rm in}\sqrt h$. From
\eqref{eq:combined-finite} and \eqref{eq:trajectory-closed} there
are constants $C,C_*>0$, uniform on the chosen tube, such that
\[
 \sup_{u\le t}\|r(u)\|_{\mathcal X}^2
 +|ha(t)-A(ht)|+|\beta(t)-B(ht)|\le C h e^{C_*ht}.
\]
Here the factor $1+ht$ has been absorbed into $e^{C_*ht}$.
Fix $0<\epsilon<1/2$ and choose
$0<\kappa<2\epsilon/C_*$. Then
$h e^{C_*ht}\le h^{1-C_*\kappa}$ for
$ht\le\kappa\log(1/h)$. For sufficiently small $h$ this gives
$\|r(t)\|_{\mathcal X}\le C h^{1/2-\epsilon}$ and keeps the
remainder strictly inside the bootstrap tube. Moreover the parameter
discrepancy is smaller than half the distance from $K$ to
$\mathbb R\setminus\operatorname{int}I$.

Start the bootstrap by local existence and the modulation construction.
At any first prospective breakdown before
$h^{-1}\min\{S_K,\kappa\log(1/h)\}$, the reduced parameter belongs
to $K$, so the fitted parameter stays strictly inside $I$. Sobolev
embedding, the smallness of the remainder, and
$1-\hat\eta_\beta\ge\min_{b\in I}\sin^2 b>0$ supply a uniform
lower bound for $1-\eta$. The full $\mathcal X$ norm is bounded,
since translates of the profiles on $I$ are uniformly bounded there.
The continuation criterion for the perturbed equation and local
uniqueness of the fit thus extend the decomposition. This excludes
breakdown before the stated time, and continuity gives the estimates
at the endpoint, including the case $B(S_K)\in\partial K$.
If the initial hypothesis is phrased as proximity to a reference
profile, the uniform modulation lemma gives $d_{\rm in}=O(\sqrt h)$;
the fitted initial parameters are still used to initialize the ODE.
\end{proof}

\begin{corollary}[Laboratory velocity and a transverse reversal]\label{cor:reversal}
Under the assumptions of Theorem~\ref{thm:main}, fix $S>0$ with
$S/h\le T_h$, and put $X_h(s)=h a(s/h)$ for $0\le s\le S$. Then
\[
 \|X_h'-A'\|_{L^\infty(0,S)}\le C_S\sqrt h,
 \qquad \|X_h'-A'\|_{L^2(0,S)}\le C_S h .
\]
In particular, take a fixed reduced orbit with
$A'(s_*)=0$ and $A''(s_*)\ne0$ at $s_*\in(0,S)$.
The fitted laboratory velocity $\dot a(s/h)=X_h'(s)$ takes opposite
signs on the two sides of $s_*$ within $C_Sh^{2/3}$ of $s_*$ in slow
time. The fitted center has a local extremum at an intervening time $s_h$,
with
\[
 |s_h-s_*|\le C_Sh^{2/3},\qquad
 |h a(s_h/h)-A(s_*)|\le C_S h .
\]
In a fixed sufficiently small neighborhood of $s_*$, every zero of $X_h'$ lies within $C_S\sqrt h$ of $s_*$. The constants for the reversal can depend on the orbit and $s_*$.
\end{corollary}

\begin{proof}
Since $X_h'(s)=\dot a(s/h)$, Proposition~\ref{prop:mod} and
\eqref{intro:trajectory} give
\[
 |X_h'(s)-A'(s)|\le C_S h+C D_R(s/h)^{1/2}.
\]
On this interval \eqref{intro:spacetime} bounds
$\sup D_R\le C_S h$ and
$\int_0^S D_R(s/h)\,ds=h\int_0^{S/h}D_R(t)\,dt\le C_Sh^2$.
These yield the two velocity estimates. At a zero of $X_h'$ near $s_*$, the uniform estimate and $|A'(s)|\ge c|s-s_*|$ give the asserted $O_S(\sqrt h)$ localization.

Near $s_*$, one has $|A'(s)|\ge c|s-s_*|$, with opposite signs
on the two sides. If $X_h'$ never takes the sign of $A'$ on a
one-sided interval of length $d$ adjacent to $s_*$, the squared
$L^2$ error there is at least $c^2d^3/3$. Taking $d=K h^{2/3}$
for a sufficiently large fixed $K$ therefore gives a point of the
corresponding strict sign on each side. Between these points the
fitted center attains an interior maximum if $A''(s_*)<0$ and an
interior minimum if $A''(s_*)>0$. At this extremum the position
estimate follows from \eqref{intro:trajectory} and
$|A(s_h)-A(s_*)|\le C_S|s_h-s_*|^2$.
\end{proof}

\section{Terminal lower bounds from the balance laws}
\label{sec:sharpness}

The expanding shelves found by Burtsev and Camassa can carry order-one area on the inverse perturbation scale \cite{BurtsevCamassa1997}. Here the localized flux balances give a lower bound in the hydrodynamic energy norm, using the change in the two core integrals.

\begin{proposition}\label{prop:sharpness}
Fix $S>0$. Suppose $\beta(t)\in I$ on $0\le t\le S/h$,
the scaled displacements $h|a(t)-a(0)|$ have a common bound $M$,
and $d=\sup_{0\le t\le S/h}\|r(t)\|_{\mathcal X}$ is bounded by
a fixed constant. Then for some $c,C_S>0$ independent of $h$,
\begin{equation}\label{eq:sharpness}
 |2\cos\beta(S/h)-2\cos\beta(0)|
 \le C_S h^{-1/2}d+C_S e^{-c/h}.
\end{equation}
\end{proposition}

\begin{proof}
Translate the spatial origin to $a(0)$ and write $V(x)$ for the translated field. All ensuing bounds are uniform under this translation. Choose $\chi\in C_c^\infty(\R)$ equal to one on a fixed neighborhood of $[-M,M]$. The continuity equation
gives
\[
 \frac{d}{dt}\int\chi(hx)\eta(t,x)\,dx
 =\int h\chi'(hx)\big(2(1-\eta)\nu+V\eta\big)\,dx .
\]
All spatial tests here are compactly supported. In particular no
integrability of the full density deficit is assumed.
The endpoint integrals differ from $2\cos\beta(t)$ by at most
$C h^{-1/2}d+C e^{-c/h}$. Indeed,
$\|\chi(h\cdot)\|_2=h^{-1/2}\|\chi\|_2$ controls the remainder
while $\int\hat\eta_\beta=2\cos\beta$ and the profile decays
exponentially, uniformly for $\beta\in I$.

Write $G_h(U)=2(1-\eta)\nu+V\eta$ for the density flux and
$Q=Q_{\beta(t)}(\cdot-a(t))$. Subtraction gives
\[
 G_h(U)-G_h(Q)
 =2(1-\eta)(\nu-Q_2)+(V-2Q_2)(\eta-Q_1).
\]
The coefficients on the right are uniformly bounded by the
$H^1$ bound for $\eta$, compactness of $I$, and boundedness of $W$.
The flux difference therefore has $L^2$ norm at most $C d$.
The pure profile flux is exponentially small on the support of
$\chi'(hx)$, whose distance to the center is at least $c_0/h$.
Since $\|h\chi'(h\cdot)\|_2=h^{1/2}\|\chi'\|_2$ and
$\|h\chi'(h\cdot)\|_1=\|\chi'\|_1$, integration over $[0,S/h]$
costs at most $C_S h^{-1/2}d+C_S e^{-c/h}$, with $c$ decreased
to absorb the time factor in the exponential tail.
This proves \eqref{eq:sharpness}.
\end{proof}

\begin{corollary}[Terminal radiation in both components]\label{cor:annular-cutoff}
Fix $S>0$ and a compact positive gray interval $I$. Consider fitted
solutions on $[0,S/h]$ with $\beta(t)\in I$ and, uniformly as
$h\downarrow0$,
\begin{align*}
 \|r(0)\|_{\X}&=o(\sqrt h),&
 \sup_{0\le t\le S/h}\|r(t)\|_{\X}&\le C_S\sqrt h,\\
 h\sup_{0\le t\le S/h}|a(t)-a(0)|&\le M_S,&
 |\beta(S/h)-\beta(0)|&\ge d_*>0.
\end{align*}
There are constants $0<\ell<L$, $c_S>0$ independent of $h$ such
that for all sufficiently small $h$ and each
$f\in\{\hat\theta(\cdot,S/h),\hat\xi(\cdot,S/h)\}$,
\begin{equation}\label{eq:terminal-annulus}
 \|f\|_{L^2(\{\ell/h\le |y|\le L/h\})}\ge c_S\sqrt h .
\end{equation}
For exact soliton initial data, Theorem~\ref{thm:main} supplies
the size and displacement bounds on each fixed slow interval it
covers. If $|B(S)-B(0)|$ has a fixed positive lower bound,
\eqref{intro:trajectory} supplies the last hypothesis.
\end{corollary}

\begin{proof}
Write $U=(\eta,\nu)$, $Q=Q_{\beta(t)}(\cdot-a(t))$ in the laboratory
frame, and let $F_1,F_2$ be the two fluxes in \eqref{eq:flow-hydro}.
The two profile components have integrals
\[
 m_1(b)=\int\hat\eta_b=2\cos b,\qquad
 m_2(b)=\int\hat\nu_b=\pi-2b .
\]
The second identity follows by integrating the negative phase
gradient of $q_b$. Each integral changes by at least
$\delta=2d_*\min_{b\in I}\sin b>0$ between the two endpoints.
The first flux difference satisfies
$\|F_1(U)-F_1(Q)\|_2\le C_S\sqrt h$ as in
Proposition~\ref{prop:sharpness}. For the second, write
$g=\log\sqrt{1-\eta}$ and $g_Q=\log\sqrt{1-Q_1}$. The gray
margin and the $H^1$ bound on the density give
$\|g-g_Q\|_{H^1}\le C\|\eta-Q_1\|_{H^1}$ by the chain rule.
In the quantum pressure $g_{xx}+g_x^2$ the second derivative
difference is $O_S(\sqrt h)$ in $H^{-1}$, and the difference
of squares is $O_S(\sqrt h)$ in $L^1$. The velocity squares
obey the same $L^1$ estimate, while $2(\eta-Q_1)+V(\nu-Q_2)$
is $O_S(\sqrt h)$ in $L^2$. Since $L^1(\R)\subset H^{-1}(\R)$,
\[
 \|F_2(U)-F_2(Q)\|_{H^{-1}}\le C_S\sqrt h.
\]

Choose $0\le\chi\le1$ smooth, equal to one on $[-1,1]$ and
supported in $[-2,2]$. For $K\ge2(M_S+1)$ put
$\chi_{h,K}(x)=\chi(h(x-a(0))/K)$. Its derivative has
$H^1$ norm at most $C\sqrt{h/K}$ for $h\le1$. Every point in
the support of this derivative lies at least $K/(2h)$ from $a(t)$.
The traveling-wave equations give $F_j(Q)=(2\sin\beta+V)Q_j$,
so integrating the two tested balance laws over $[0,S/h]$ gives
\[
 \left|\int\chi_{h,K}(x)[U_j(S/h,x)-U_j(0,x)]\,dx\right|
 \le C_S K^{-1/2}+C_S e^{-cK/h},\qquad j=1,2 .
\]
Here the first term bounds the flux difference. The exponential
term uses the uniform profile decay and
$\|\chi_{h,K}'\|_1=\|\chi'\|_1$, after decreasing $c$ to absorb
the time factor. The tested profile integrals at the endpoints
are $m_j(\beta(t))+O(e^{-cK/h})$. The tested initial remainder is
$o(1)$ for each fixed $K$, since $\|\chi_{h,K}\|_2=O(\sqrt{K/h})$.
Choose $K$ large enough that $C_SK^{-1/2}<\delta/4$, then $h$
small enough to obtain
\[
 \left|\int\chi_{h,K}(a(S/h)+y)\,r_j(y,S/h)\,dy\right|
 \ge\delta/2 .
\]
The part with $|y|<\ell/h$ is at most $C_S\sqrt{2\ell}$ by
Cauchy--Schwarz. Choose $\ell$ so that this is at most $\delta/4$.
The support of the remaining test lies in $|y|\le(2K+M_S)/h$.
A second application of Cauchy--Schwarz with $L=2K+M_S$
proves \eqref{eq:terminal-annulus} for each component.
\end{proof}

\begin{corollary}[Sharpness for exact soliton data]\label{cor:matching}
Take $W(A)=\sech^2 A$ and $U_0=Q_{\pi/4}(\cdot+1/h)$.
For every fixed $S>0$ the terminal field obeys
\[
 c_S\sqrt h\le \|r(\cdot,S/h)\|_{\X}
 \le \sup_{0\le t\le S/h}\|r(\cdot,t)\|_{\X}\le C_S\sqrt h
\]
for all sufficiently small $h$. The cutoff estimate
\eqref{eq:terminal-annulus} applies directly at every $S>0$
except the unique return time $S_{\rm ret}$ with $A(S_{\rm ret})=1$.
\end{corollary}

\begin{proof}
Here $A(0)=-1$, $B(0)=\pi/4$, and the energy level
$e=\sqrt2/3+(\pi/4-1/2)\sech^2(1)$ lies in $(0,4/3)$.
The aiding-field orbit in Section~1 stays in a compact gray interval,
and $A'$ is bounded below by a positive constant. Hence it crosses
$A=1$ once. Since $P(\pi/4)>0$ and the energy at a fixed $A$
is strictly decreasing in $B$, the angle returns to $\pi/4$
precisely when $W(A)=W(-1)$, namely at $A=1$ after the initial time.
Theorem~\ref{thm:main} gives the upper bound on each fixed interval,
and Corollary~\ref{cor:annular-cutoff} gives the terminal lower bound
for $S\ne S_{\rm ret}$. At $S_{\rm ret}$ the endpoint profile
integrals agree. Since $W'(-1)>0$, \eqref{ode} gives $B'(0)>0$.
Corollary~\ref{cor:acoustic-lower}, proved from Theorem~\ref{thm:acoustic} below, gives the annular lower bound at the return time as well.
\end{proof}

\section{The macroscopic radiation profile}\label{sec:acoustic}

Perturbative analyses have described shelves forming on both sides of a dark soliton. Burtsev and Camassa found an order-one integrated shelf on the inverse-perturbation scale for weak damping or gain \cite{BurtsevCamassa1997}. Ablowitz, Nixon, Horikis and Frantzeskakis used conservation laws and matched inner and outer expansions to compute shelf propagation \cite{AblowitzNixonHorikisFrantzeskakis2011}. Ossi, Prinari and Yang derived shelf and parameter formulas using squared eigenfunctions \cite{OssiPrinariYang2026}. A comparison of spectral and perturbed conservation-law approaches appears in \cite{OssiPrinariHorikisFrantzeskakis2026}. For the prescribed field in \eqref{intro:forced}, the energy and local virial estimates justify the macroscopic limit directly for nonlinear energy solutions.

Fix $A_0\in\R$, $B_0\in(0,\pi/2)$, and let $(A,B)$ solve
\eqref{ode} with these initial values. Fix $S>0$ such that
$B([0,S])$ lies in a compact subinterval of the positive gray branch.
We use laboratory coordinates at both scales. For $s=ht$, $z=hx$, set
\begin{equation}\label{eq:acoustic-scale}
 q_h(s,x)=Q_{\beta(s/h)}(x-a(s/h)),\qquad
 d_h(s,x)=U(s/h,x)-q_h(s,x),\qquad
 R_h(s,z)=h^{-1}d_h(s,z/h).
\end{equation}
Thus $d_h$ is the translate of the fitted remainder in
\eqref{intro:orth}. The integrated density deficit and negative phase
gradient of the reduced profile form the vector
\begin{equation}\label{eq:acoustic-masses}
 m(s)=\int_{\R} Q_{B(s)}(y)\,dy
     =\big(2\cos B(s),\,\pi-2B(s)\big)^T .
\end{equation}
The first component follows from \eqref{intro:Q}, and the second from the phase difference of $q_B$ across the line. These integrals concern the core. For the full energy-class solution we will test both local conservation laws against compactly supported functions. The renormalized momentum $P(B)$ enters the force law \eqref{ode}.

\begin{theorem}[Acoustic radiation from exact soliton data]\label{thm:acoustic}
Suppose $U(0,x)=Q_{B_0}(x-A_0/h)$. Let $R$ solve
\begin{equation}\label{eq:acoustic-system}
 \partial_s R+\partial_z(C(z)R)=-m'(s)\delta_{A(s)},\qquad
 C(z)=W(z)\operatorname{Id}+2\sigma,\qquad R(0)=0.
\end{equation}
This equation has a unique distributional solution in
$C([0,S];L^2(\R)^2)$ with the indicated trace. For all sufficiently
small $h$,
\begin{equation}\label{eq:acoustic-rate}
 \sup_{0\le s\le S}\|R_h(s)-R(s)\|_{H^{-3}(\R)^2}\le C_S h.
\end{equation}
For each $g\in L^2(\R)^2$ the pairings
$\langle R_h(s)-R(s),g\rangle$ tend to zero uniformly in $s\in[0,S]$.

More explicitly, for $\epsilon\in\{-1,1\}$ set
$R^\epsilon=R_1+\epsilon R_2$, $m^\epsilon=m_1+\epsilon m_2$.
We write $R^{+},R^{-}$ for the choices $\epsilon=+1,-1$, and likewise for $m$. Let $\Phi_\epsilon(s;\tau,z)$ be the position at time $s$ of the
characteristic starting from $z$ at time $\tau$ with velocity
$W+2\epsilon$. Write
\[
 Z_\epsilon(s,\tau)=\Phi_\epsilon(s;\tau,A(\tau)),\qquad
 J_\epsilon(s,\tau)=\partial_z\Phi_\epsilon(s;\tau,A(\tau)).
\]
At each $s>0$, $R^\epsilon(s,\cdot)$ is supported on the interval
with endpoints $A(s)$ and $\Phi_\epsilon(s;0,A_0)$. Its density in the
interior is given by
\begin{equation}\label{eq:acoustic-profile}
 R^\epsilon(s,Z_\epsilon(s,\tau))
 =-\frac{(m^\epsilon)'(\tau)}
        {J_\epsilon(s,\tau)\,|2\sin B(\tau)-2\epsilon|},
 \qquad 0<\tau<s .
\end{equation}
\end{theorem}

\begin{proof}
Theorem~\ref{thm:main} applies on $[0,S/h]$ for small $h$. Its
estimates \eqref{intro:spacetime}--\eqref{intro:trajectory} give
\begin{align}
 \sup_{s\le S}\|d_h(s)\|_{\X}&\le C_S\sqrt h,&
 \int_0^S D_R(\tau/h)\,d\tau&\le C_Sh^2,\notag\\
 \sup_{s\le S}\big(|ha(s/h)-A(s)|+|\beta(s/h)-B(s)|\big)&\le C_Sh .
 \label{eq:acoustic-input}
\end{align}
In particular $\sup_s\|R_h(s)\|_2\le C_S$. We keep the initial
radiation term in the proof for use in Corollary~\ref{cor:acoustic-prepared}.

\emph{The concentrated profile and its flux.}
The components of $Q_b$ have bounded first absolute moments,
uniformly on the fitted gray interval. Translating to $y=x-a(s/h)$
in the pairing against $g(hx)$ and using \eqref{eq:acoustic-input}
gives, componentwise,
\[
 \left|\int q_h(s,x)g(hx)\,dx-m(s)g(A(s))\right|
 \le C_Sh\|g\|_{W^{1,\infty}} .
\]
Let $F_h(U)$ denote the flux in \eqref{eq:flow-hydro}, so that
$U_t+\partial_xF_h(U)=0$, and write $F_{\mathrm{free}}$ for the unforced
flux. A traveling profile satisfies $F_{\mathrm{free}}(Q_b)=c(b)Q_b$ by integration
of its equation and decay at infinity. Thus
$F_h(q_h)=(c(\beta(s/h))+W(hx))q_h$. The same pairing estimate
holds with $q_h,m$ replaced by $F_h(q_h),A'm$, since
$A'=c(B)+W(A)$ and $W,W'$ are bounded.

\emph{The remaining flux.}
For this computation put $q=q_h(s,\cdot)$, $d=d_h(s,\cdot)$ and
$\rho=\sqrt{1-\eta(s/h,\cdot)}$. Set
\[
 l=\rho_x/\rho=-\eta_x/[2(1-\eta)],\qquad
 l_q=-q_{1,x}/[2(1-q_1)],\qquad k=l-l_q .
\]
The identity $\rho_{xx}/\rho=\partial_x l+l^2$ holds in distributions
by the $H^1$ multiplier rule and the uniform lower bound on $\rho$.
Subtracting $F_h(q)+C(hx)d$ from $F_h(q+d)$ leaves
\[
 K_1=-2(q_1d_2+q_2d_1+d_1d_2),\qquad
 K_2=\partial_x k+2l_q k+k^2-2q_2d_2-d_2^2 .
\]
In particular $|k|\le C(|d_{1,x}|+|q_{1,x}d_1|)$, so
$\|k\|_2\le C_S\sqrt h$. Products of two radiation factors have
$L^1_x$ norm $O_S(h)$. Each mixed product containing $q_1,q_2$ or
$l_q$ has $L^1_x$ norm at most $C\sqrt{D_R(s/h)}$ by exponential
decay and weighted Cauchy--Schwarz. The distributional derivative
satisfies
$|\langle\partial_x k,g(h\cdot)\rangle|
 \le\sqrt h\,\|k\|_2\|g'\|_2$.
It follows from \eqref{eq:acoustic-input} that for vector tests $g$
and $0\le s\le S$,
\begin{equation}\label{eq:acoustic-flux-defect}
 \int_0^s\!|\langle K(\tau),g(\tau,h\cdot)\rangle|\,d\tau
 \le C_S h\sup_{\tau\le s}
   \big(\|g(\tau)\|_\infty+\|\partial_z g(\tau)\|_2\big).
\end{equation}
Indeed $\int_0^s\sqrt{D_R(\tau/h)}\,d\tau\le C_Sh$.

\emph{Identification of the limit.}
Changing variables in \eqref{eq:flow-hydro} gives the two
conservation laws in $(s,z)$, with the total field and flux both
multiplied by $h^{-1}$. Test these laws between $0$ and $s$ with a smooth vector
$\phi(\tau,z)$ compactly supported in $z$. Use the two concentrated
profile estimates on $\phi,\phi_\tau,\phi_z$ and
\eqref{eq:acoustic-flux-defect} on $g=\phi_z$. Integrating the
derivative of $m(\tau)\cdot\phi(\tau,A(\tau))$ gives
\begin{align}
 \langle R_h(s),\phi(s)\rangle-\langle R_h(0),\phi(0)\rangle
 ={}&\int_0^s\!\langle R_h,\phi_\tau+C\phi_z\rangle\,d\tau\notag\\
 &-\int_0^s\!m'(\tau)\cdot\phi(\tau,A(\tau))\,d\tau
       +\mathcal E_h .\label{eq:acoustic-weak}
\end{align}
Here $|\mathcal E_h|\le C_Sh$ times the supremum of the
$W^{1,\infty}$ norms of $\phi,\phi_\tau,\phi_z$ and the $L^2$
norm of $\phi_{zz}$. The weak law uses only spatial distributions
for the physical flux; its time endpoints are justified by the
continuity from Proposition~\ref{prop:flow}.

For $f\in C_c^\infty(\R)$ choose
$\phi(\tau,z)=(1,\epsilon)^T f(\Phi_\epsilon(s;\tau,z))$.
Then $\phi_\tau+C\phi_z=0$. The associated error is bounded by
$C_S h\|f\|_{H^3}$ using the bounded derivatives and Jacobians of
the characteristic flow and $H^3(\R)\hookrightarrow W^{2,\infty}(\R)$.
For the exact data $R_h(0)=0$. Thus the candidate $R^\epsilon$ is
specified by
\begin{equation}\label{eq:acoustic-duality}
 \langle R^\epsilon(s),f\rangle
    =-\int_0^s (m^\epsilon)'(\tau)f(Z_\epsilon(s,\tau))\,d\tau .
\end{equation}
Differentiation against smooth tests verifies its scalar
conservation law and initial trace in \eqref{eq:acoustic-system};
the same backward tests give uniqueness.

To evaluate \eqref{eq:acoustic-duality}, differentiate the
characteristic with respect to its starting time and point. This gives
\[
 \partial_\tau Z_\epsilon(s,\tau)
   =J_\epsilon(s,\tau)(2\sin B(\tau)-2\epsilon),\qquad
 J_\epsilon(s,\tau)
   =\exp\!\left(\int_\tau^s W'(\Phi_\epsilon(u;\tau,A(\tau)))\,du\right).
\]
The absolute derivative is bounded above and away from zero on
$0\le\tau\le s\le S$. Changing variables yields
\eqref{eq:acoustic-profile}, with the $+1$ characteristic interval to the
right of $A(s)$ and the $-1$ interval to the left. In particular the
densities are bounded, with supports in a fixed compact interval.
Their endpoints depend continuously on $s$, and the formula in the
interior depends continuously on $s$ by inversion of the strictly
monotone map. Dominated convergence gives continuity in $L^2$,
also at $s=0$ where both intervals have zero length.
Duality and the $H^3$ bound for the error prove
\eqref{eq:acoustic-rate}. The uniform $L^2$ bounds on $R_h,R$ and
density of smooth tests in $L^2$ then give the asserted weak
convergence.
\end{proof}

\begin{corollary}[Terminal radiation in each component]\label{cor:acoustic-lower}
Under the hypotheses of Theorem~\ref{thm:acoustic}, for $j=1,2$,
\begin{equation}\label{eq:acoustic-lower}
 \liminf_{h\downarrow0}h^{-1/2}\|d_{h,j}(S)\|_{L^2_x}
 \ge \frac12\left(
   \sum_{\epsilon=\pm1}\int_0^S
   \frac{|(m^\epsilon)'(\tau)|^2}
        {J_\epsilon(S,\tau)|2\sin B(\tau)-2\epsilon|}\,d\tau
 \right)^{1/2} .
\end{equation}
If $B$ is nonconstant on $[0,S]$, this coefficient is positive, even
when $B(S)=B(0)$. In that case there exist $0<b_1<b_2$ and $k_S>0$
such that, for all sufficiently small $h$, both components satisfy
\begin{equation}\label{eq:acoustic-annular}
 \|d_{h,j}(S)\|_{L^2(\{x:\ b_1/h<|x-a(S/h)|<b_2/h\})}
     \ge k_S\sqrt h .
\end{equation}
\end{corollary}
\begin{proof}
Weak lower semicontinuity gives a lower bound by $\|R_j(S)\|_2$
in \eqref{eq:acoustic-lower}. The two characteristic intervals
have disjoint interiors. Hence
$\|R_j(S)\|_2^2=(\|R^{+}(S)\|_2^2+\|R^{-}(S)\|_2^2)/4$ for either $j$. Changing variables in
\eqref{eq:acoustic-profile} gives the claimed coefficient.
Since $(m^\epsilon)'=-2(\sin B+\epsilon)B'$, a nonconstant $B$
has a closed subinterval of $(0,S)$ of positive length where both
sources are bounded away from zero. Its image under either
$Z_\epsilon(S,\cdot)$ is compactly contained in the corresponding
characteristic interval and separated from $A(S)$. Both $R_j(S)$
have positive $L^2$ norm on that image. Apply weak lower
semicontinuity there and choose $b_1,b_2$ with a positive distance
margin at both endpoints. Since $ha(S/h)=A(S)+O_S(h)$, the
corresponding physical region lies inside the annulus in
\eqref{eq:acoustic-annular} for small $h$.
\end{proof}

\begin{corollary}[Prepared initial radiation]\label{cor:acoustic-prepared}
Suppose the data have fitted initial parameters $ha(0)=A_0$,
$\beta(0)=B_0$ and $\|d_h(0)\|_{\X}\le C_0\sqrt h$.
For $R^0\in L^2(\R)^2$, let $R^{\mathrm{in}}$ be the solution of
\eqref{eq:acoustic-system} with initial value $R^0$ in place of zero.
Then
\begin{equation}\label{eq:acoustic-prepared}
 \sup_{s\le S}\|R_h(s)-R^{\mathrm{in}}(s)\|_{H^{-3}(\R)^2}
 \le C_{S,C_0}\big(h+\|R_h(0)-R^0\|_{H^{-3}(\R)^2}\big).
\end{equation}
If $R_h(0)\rightharpoonup R^0$ in $L^2$, then $R_h(s)$ converges to
$R^{\mathrm{in}}(s)$ weakly in $L^2$ against every fixed test,
uniformly on $[0,S]$. In particular the error is $O_{S,C_0}(h)$ after subtracting the homogeneous acoustic evolution of $R_h(0)$ and the zero-data emitted field. With $R^{\mathrm{em}}$ the zero-initial-data
solution of Theorem~\ref{thm:acoustic} and $(R^0)^\epsilon=R^0_1+\epsilon R^0_2$,
the limit has the component formula
\[
 (R^{\mathrm{in}})^\epsilon(s,z)=(R^{\mathrm{em}})^\epsilon(s,z)
 +\partial_z\Phi_\epsilon(0;s,z)\,(R^0)^\epsilon(\Phi_\epsilon(0;s,z)).
\]
\end{corollary}
\begin{proof}
All estimates through \eqref{eq:acoustic-weak} hold with constants
depending on $C_0$, because \eqref{intro:spacetime} still gives
\eqref{eq:acoustic-input}. The initial term against the backward
test is now $\langle R_h^\epsilon(0),f(\Phi_\epsilon(s;0,\cdot))\rangle$,
where $R_h^\epsilon=R_{h,1}+\epsilon R_{h,2}$. Replacing the
initial factor by $(R^0)^\epsilon$ costs at most
$C_S\|R_h(0)-R^0\|_{H^{-3}}\|f\|_{H^3}$. The change of variables
$z=\Phi_\epsilon(s;0,y)$ in the resulting pairing gives the stated
transported initial density and \eqref{eq:acoustic-prepared}. The constants are independent of $R^0$, so one may take $R^0=R_h(0)$ separately for each $h$. An incoming field can therefore partially cancel the emitted field in the terminal profile.
Under weak $L^2$ convergence alone, the initial pairings converge
uniformly in $s$ for each smooth compactly supported $f$, since the
family $f(\Phi_\epsilon(s;0,\cdot))$, $0\le s\le S$, is compact in
$L^2$. The uniform $L^2$ bounds extend this conclusion to every fixed
$L^2$ test.
\end{proof}

\begin{figure}[htbp]
\centering
\includegraphics[width=.93\linewidth]{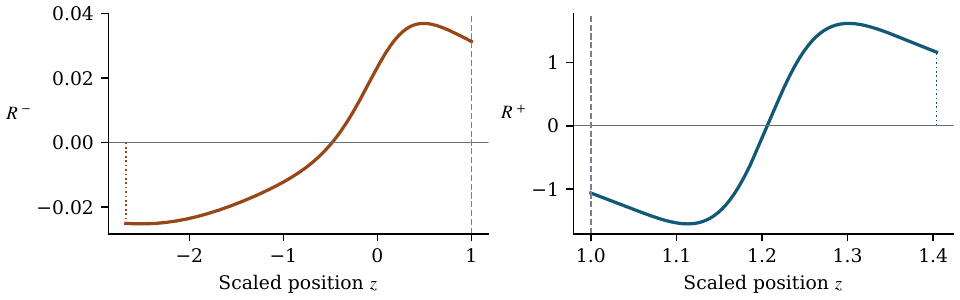}
\caption{The two limiting acoustic fields for $W(A)=\sech^2 A$, $A(0)=-1$ and $B(0)=\pi/4$, evaluated from the characteristic formula at the slow time $S\simeq0.90262$ with $A(S)=1$. The dashed line marks the final core. The panels use different vertical scales. Symmetry gives $B(S)=B(0)$. Each component $R^{\pm}$ has zero signed integral and positive $L^2$ norm.}
\label{fig:acoustic}
\end{figure}

\smallskip\noindent\emph{An excursion with zero net emission in each family.}\enspace
Figure~\ref{fig:acoustic} follows the aiding field in the first phase portrait. The incoming orbit has $A'>0$. Its angle is determined uniquely by the energy and the local field value, hence its values agree at $A=-1$ and $A=1$. The total signed integral of each emitted field is $m^{\pm}(0)-m^{\pm}(S)=0$. The contributions from the two flanks occupy different positions on the corresponding characteristic interval. Their squared norms therefore retain the effect of the excursion.

\smallskip\noindent\emph{Further questions.}\enspace
The weak limit describes radiation on the spatial scale $h^{-1}$. A strong $L^2$ limit, or an $o(\sqrt h)$ approximation after subtracting a radiation corrector, would require control of energy on shorter scales. The limiting density can jump at an acoustic front, and a corrector measured in $H^1$ would need to resolve that transition. For the opposing field, the reduced saddle suggests studying the existence and stability type of a pinned PDE profile and the possible capture of a critical incident orbit for fixed $h$. The comparison proved here applies on compact gray bands. At the black boundary, a continuation theory would also have to treat the reciprocal in \eqref{intro:complex}.

\bibliographystyle{amsalpha}
\bibliography{references}
\end{document}